\documentclass[12pt]{amsart}

\usepackage{amssymb}
\usepackage[all]{xy}
\usepackage{eucal}
\usepackage{xfrac}
\usepackage{relsize}
\usepackage{stmaryrd}

\usepackage[backend=biber,
            isbn=false,
            doi=false,
            maxbibnames=5,
            giveninits=true,
            style=alphabetic,
            citestyle=alphabetic]{biblatex}
\usepackage{url}
\DeclareFieldFormat{postnote}{\mknormrange{#1}}
\DeclareFieldFormat{volcitepages}{\mknormrange{#1}}
\DeclareFieldFormat{multipostnote}{\mknormrange{#1}}
\renewbibmacro{in:}{}
\bibliography{bibliography.bib}

\usepackage{tikz}
\usepackage{tikz-cd}
\usetikzlibrary{calc,matrix}
\tikzset{w/.style={circle, draw,inner sep=1pt},b/.style={circle,draw,fill,inner sep=2pt}, s/.style={rectangle, draw,inner sep=3pt}}

\usepackage[colorlinks=true]{hyperref}
\usepackage{cleveref}
\usepackage{textcomp}

\newcounter{conj}
\newcounter{intro}
\newtheorem{thm}{Theorem}[section]
\crefname{thm}{Theorem}{Theorem}
\newtheorem{prop}[thm]{Proposition}
\crefname{prop}{Proposition}{Proposition}

\crefname{constr}{Construction}{Construction}
\newtheorem{lm}[thm]{Lemma}
\crefname{lm}{Lemma}{Lemma}

\crefname{cor}{Corollary}{Corollary}

\crefname{conjecture}{Conjecture}{Conjecture}

\crefname{question}{Question}{Question}

\newtheorem{thmintro}[intro]{Theorem}

\theoremstyle{definition}
\newtheorem{defn}[thm]{Definition}
\crefname{defn}{Definition}{Definition}

\theoremstyle{remark}
\newtheorem{remark}[thm]{Remark}
\crefname{remark}{Remark}{Remark}
\newtheorem{example}[thm]{Example}
\crefname{example}{Example}{Example}
\newtheorem{warning}[thm]{Warning}
\crefname{warning}{Warning}{Warning}

\newcommand{\defeq}{\mathbin{:\!=}}
\newcommand{\defterm}[1]{\textbf{\emph{#1}}}

\newcommand{\eps}{\varepsilon}
\newcommand{\id}{\mathrm{id}}

\newcommand{\calF}{\mathcal{F}}
\newcommand{\calG}{\mathcal{G}}
\newcommand{\calH}{\mathcal{H}}
\newcommand{\calB}{\mathcal{B}}
\newcommand{\calC}{\mathcal{C}}
\newcommand{\calD}{\mathcal{D}}
\newcommand{\calE}{\mathcal{E}}
\newcommand{\calL}{\mathcal{L}}
\newcommand{\calS}{\mathcal{S}}
\newcommand{\calO}{\mathcal{O}}
\newcommand{\calM}{\mathcal{M}}

\newcommand{\Sp}{\mathcal{S}}
\newcommand{\bP}{\mathbb{P}}
\newcommand{\bD}{\mathbb{D}}
\newcommand{\bA}{\mathbb{A}}

\newcommand{\fg}{\mathfrak{g}}
\newcommand{\rB}{\mathrm{B}}
\newcommand{\rZ}{\mathrm{Z}}
\newcommand{\Z}{\mathbf{Z}}
\newcommand{\conv}{\mathrm{conv}}
\newcommand{\fet}{\mathrm{fet}}
\newcommand{\gr}{\mathrm{gr}}
\newcommand{\fil}{\mathrm{fil}}
\newcommand{\cfil}{\wedge\mathrm{fil}}
\renewcommand{\int}{\mathrm{int}}
\newcommand{\ext}{\mathrm{ext}}
\newcommand{\op}{\mathrm{op}}
\newcommand{\nd}{\mathrm{nd}}
\newcommand{\conn}{\mathrm{conn}}

\newcommand{\ap}{\mathrm{ap}}
\newcommand{\lafp}{\mathrm{lafp}}
\newcommand{\symp}{\mathrm{symp}}
\renewcommand{\cot}{\mathbb{L}}
\renewcommand{\tan}{\mathbb{T}}

\newcommand{\colim}{\operatorname{colim}\limits}
\newcommand{\Map}{\mathrm{Map}}
\newcommand{\MapSt}{\underline{\mathrm{Map}}}
\newcommand{\Fun}{\mathrm{Fun}}
\newcommand{\Hom}{\mathrm{Hom}}
\newcommand{\fib}{\mathrm{fib}}
\newcommand{\cofib}{\mathrm{cofib}}
\newcommand{\adic}{\mathrm{adic}}

\newcommand{\BV}{\mathrm{BV}}
\newcommand{\CE}{\mathrm{CE}}
\newcommand{\KT}{\mathrm{KT}}
\newcommand{\DR}{\mathrm{DR}}
\newcommand{\Pol}{\mathrm{Pol}}
\newcommand{\Pois}{\mathrm{Pois}}
\newcommand{\Cois}{\mathrm{Cois}}
\newcommand{\Symp}{\mathrm{Symp}}
\newcommand{\Lag}{\mathrm{Lag}}
\newcommand{\Isot}{\mathrm{Isot}}
\newcommand{\dAff}{\mathrm{dAff}}
\newcommand{\stCrit}{\mathrm{Crit}^{\mathrm{st}}}
\newcommand{\dCrit}{\mathrm{dCrit}}

\newcommand{\Vect}{\mathrm{Vect}}
\newcommand{\dPSt}{\mathrm{dPSt}}
\newcommand{\CAlg}{\mathrm{CAlg}}
\newcommand{\Alg}{\mathrm{Alg}}
\newcommand{\Mod}{\mathrm{Mod}}
\newcommand{\Fol}{\mathrm{Fol}}
\newcommand{\Spec}{\mathrm{Spec}}
\newcommand{\Sym}{\mathrm{Sym}}
\newcommand{\Thick}{\mathrm{Thick}}
\newcommand{\Lie}{\mathrm{Lie}}
\newcommand{\Com}{\mathrm{Com}}

\newcommand{\QCoh}{\mathrm{QCoh}}

\newcommand{\red}{\mathrm{red}}
\newcommand{\dR}{\mathrm{dR}}

\newcommand{\vari}{\textendash}

\begin{document}
\title{AKSZ construction for shifted Poisson structures}
\author{Nikola Tomi\'c}
\begin{abstract}
	We prove the AKSZ theorem for shifted Poisson structures: if $X$ is an $n$-shifted Poisson derived stack, and $Y$ a $d$-oriented derived stack, then the mapping stack
	\[\MapSt(Y,X)\]
	is naturally endowed with an $(n-d)$-shifted Poisson structure. For this, we prove that the data of an $n$-shifted Poisson structure on a derived Artin stack is equivalent to the data of an $(n+1)$-shifted Lagrangian thickening of it. We also extend the definition of shifted Poisson structures to derived prestacks having a deformation theory and give two applications, one for mapping stacks with a non-proper source and one in BV formalism.
\end{abstract}
\maketitle
\tableofcontents
\section*{Introduction}
\subsection*{Background}
\par Poisson geometry is the study of Poisson manifolds, that is, manifolds $M$ given with a biderivation
\[\{\vari,\vari\} : C^\infty(M)\otimes C^\infty(M) \rightarrow C^\infty(M)\]
which is also a Lie bracket. The study of such objects goes back to Lie, Cartan, and Dirac. They were motivated by the Hamiltonian formalism in physics, where a phase space always comes with a Poisson structure. A celebrated result in Poisson geometry claims that, for a Poisson manifold $M$, there exists a deformation quantization on its ring of functions $C^\infty(M)$, that is, an associative structure on $C^\infty(M)[\![\hbar]\!]$ that specializes to the usual commutative structure when $\hbar = 0$. This theorem was proven in the nineties by Kontsevich (\cite{defquant}). Nowadays, the tools used to prove this theorem have become an active subject of research.
\par In algebraic geometry, it is usual to study the structure of moduli stacks. For instance, the moduli stack of genus $0$ curves is related to physics (it can be seen as a stack of trajectories of a particle in the space-time), and to Galois theory (where it is possible to obtain nice representations of the absolute Galois group of $\mathbb{Q}$). It is usually really difficult to study those objects because they have a lot of singularities.
\par In order to deal with singularities, a powerful theory wes created: derived geometry. Derived geometry was first developed by To\"en--Vezzosi, Lurie, and Gaitsgory--Rozenblyum (\cite{HAGII,SAG,GRI,GRII}) in the Noughties, it allows to manipulate singular objects as if they had no singularities. Sometimes, this phenomenon is referred to as ``Kontsevich's hidden smoothness." So, it is possible to prove that all usual moduli stacks come with a natural derived structure. Once this is settled, we can use tools from derived geometry to study them. One of the great successes of derived geometry is the theory of shifted symplectic structures that have been developed by Pantev--To\"en--Vaqui\'e--Vezzossi (\cite{PTVV}) and has found many applications in enumerative geometry (see \cite{DTsheaf, HHRII, symppull}) and in representation theory (see \cite{akszsymp}). A derived stack $X$ is said to have an $n$-shifted symplectic structure if it is given with a $2$-form
\[\omega_2 \in \wedge^2\cot_X[n]\]
that is non-degenerate: the natural morphism
\[\Theta_{\omega_2} : \tan_X \rightarrow \cot_X[n]\]
is an equivalence. This $2$-form is also assumed to be closed in a higher categorical sense: it is given with extra forms
\[\omega_i \in \wedge^i\cot_X[n-i]\]
such that
\[\rm{d}_\dR \omega_i = \rm{d}\omega_{i+1}.\]
This definition is really convenient because it allows a flexible approach to the theory of symplectic forms in the realm of algebraic geometry. Indeed, a lot of functors usually are not exact, but if we derive everything and view them in the scope of higher category theory, they become exact and a lot of constructions can be done (here we don't assume that $\omega_2$ is strictly closed, but as a counterpart we have to deal with higher homotopical data; it turns out that this data is really useful to show the existence of new shifted symplectic structures). It usually makes sense to talk about shifted symplectic structures for $X$ assumed to be an Artin stack locally of finite presentation, but it was recently extended by Calaque--Safronov (\cite{sympgrpd}) to prestacks that have a perfect cotangent complex. In this setting, the AKSZ construction (named after Alexandrov–Kontsevich–Schwarz–Zaboronsky \cite{AKSZ}) refers to a process that, given a derived stack $Y$ having a nice $d$-dimensional duality theory, and a derived stack $X$ with an $n$-shifted symplectic form $\omega$, associates an $(n-d)$-shifted symplectic form on the mapping stack
\[\MapSt(Y,X).\]
When $Y$ is the Betti stack of a $d$-dimensional connected compact oriented manifold $M$, and $X$ is the classifying stack $\rB G$ of a reductive group $G$, we get a $(2-d)$-shifted symplectic structure on a derived enhancement of the moduli stack of $G$-local systems on $M$
\[\mathbb{R}\rm{Loc}_G(M).\]
When $M$ is a surface, this is a slick way to recover the Atiyah--Bott--Goldman symplectic structure on the character stack of $M$.
\par Shifted Poisson geometry has been developed by Calaque--Pantev--To\"en--Vaqui\'e--Vezzosi (\cite{CPTVV}) and is a generalization of classical Poisson geometry in the derived setting. A $n$-shifted Poisson structure on a derived Artin stack $X$ locally of finite presentation is the data of a biderivation
\[\{\vari,\vari\} : \calO(X) \otimes \calO(X) \rightarrow \calO(X)[-n]\]
and higher operations witnessing the fact that $\{\vari,\vari\}$ is a Lie bracket up to homotopy. Classically, a symplectic manifold is a special case of a Poisson manifold. Namely, the datum of a non-degenerate Poisson brackets on a manifold is the same datum as a symplectic form on it. There is the analog statement in derived geometry:
\[\Pois^\nd(X,n) \simeq \Symp(X,n)\]
where on the left we have the space of non-degenerate $n$-shifted Poisson structures on $X$, and on the right the space of $n$-shifted symplectic structures on $X$. This theorem is highly non-trivial in this setting because one has to deal with higher homotopical data; it is one of the main theorem of \cite{CPTVV} and, independently, of \cite{poisnd}. Another important theorem of \cite{CPTVV} is the existence of deformation quantizations of the category
\[\rm{Perf}(X)\]
of perfect sheaves on $X$, when $X$ has an $n$-shifted Poisson structure (for $n\neq 0$). Thanks to the AKSZ construction for shifted symplectic structures, we get many interesting quantizations of moduli stacks that way. It happens sometimes that some moduli stacks don't naturally have a shifted symplectic structure but a shifted Poisson structure, it is then natural to ask if there is an analog of the AKSZ construction for shifted Poisson structures. Using the tools of formal derived geometry of Gaitsgory--Rozenblyum (\cite{GRI, GRII}), and formal integration of derived foliations of Brantner--Magidson--Nuiten and Fu (\cite{deffol, follie}), the author proved in \cite{lagthick}, a restricted version of the AKSZ construction, where $X$ is a derived scheme, and $Y = M_B$ is the Betti stack of a manifold $M$. In this paper, we establish the general AKSZ construction:
\begin{thmintro}[\ref{thm:aksz}]\label{thm:introaksz}
	If $Y$ is an $\calO$-compact $d$-oriented prestack and $X$ an $n$-shifted Poisson derived Artin stack such that $\MapSt(Y,X)$ is locally of finite presentation. Then the mapping prestack
	\[\MapSt(Y,X)\]
	has a natural $(n-d)$-shifted Poisson structure.
\end{thmintro}
\subsection*{Strategy}
\par The proof follows a similar strategy as in \cite{lagthick}. For an $n$-shifted Poisson derived scheme, it was established that it is possible to construct a formal thickening
\[X \rightarrow X^\symp\]
with an $(n+1)$-shifted Lagrangian structure on it. Conversely, such a thickening provides an $n$-shifted Poisson structure on $X$. The idea is then to apply $\MapSt(Y,\vari)$ on this morphism, and use an AKSZ construction for Lagrangian morphisms (\cite[2.35. (2)]{sympgrpd}) to get an $(n-d+1)$-shifted Lagrangian thickening of $\MapSt(Y,X)$. Unfortunately, the theorem proved in \cite{lagthick} was too restrictive, because in really few cases $\MapSt(Y,X)$ is still a derived scheme. First, we extend the definition of shifted Poisson structures on $X$ for $X$ being a prestack only assumed to have a deformation theory and be locally of finite presentation, and then we extend the theorem for such an $X$. Thus proving a conjecture of Calaque (\cite[Conjecture 3.4.]{cotstack}):
\begin{thmintro}[\ref{thm:lagthick}]\label{thm:intro}
	If $X$ has a deformation theory and is locally of finite presentation, then there is an equivalence:
	\[\Pois(X,n) \simeq \Lag\Thick(X,n+1)\]
	between the spaces of $n$-shifted Poisson structures on $X$ and $(n+1)$-shifted Lagrangian thickenings of $X$.
\end{thmintro}
\par The AKSZ construction is then almost immediate: If $X$ has a deformation theory, is locally of finite presentation, and has an $n$-shifted Poisson structure, we extract an $(n+1)$-shifted Lagrangian thickening
\[X \rightarrow X^\symp\]
and apply $\MapSt(Y,\vari)$ on it:
\[f : \MapSt(Y,X) \rightarrow \MapSt(Y,X^\symp).\]
We get an $(n-d+1)$-shifted Lagrangian on it. Now, it is automatic to check that $\MapSt(Y,X)$ has a deformation theory. By assumption, $\MapSt(Y,X)$ is locally of finite presentation, so we can apply theorem \ref{thm:intro} in order to extract a shifted Poisson structure on $\MapSt(Y,X)$. The map $f$ may be not a formal thickening, however, we can take its formal completion. That way we get an $(n-d+1)$-shifted Lagrangian thickening of $\MapSt(Y,X)$ and so an $(n-d)$-shifted Poisson structure on it. The hard part actually lies in the proof of theorem \ref{thm:intro}.
\par Shifted Poisson structures on an arbitrary prestack $X$ with a cotangent complex are defined using formal localization, we view $X$ as a family over its de Rham prestack $\pi : X \rightarrow X_\dR$. For $\Spec(A)\rightarrow X_\dR$, we set $X_A$ to be the pullback of $\Spec(A)$ along $\pi$. Note that the map $\Spec(A) \rightarrow X_\dR$ induces a map $\Spec(A^\red) \rightarrow X$, so we can picture those constructions as follows:
\[\begin{tikzcd}
	\Spec(A^\red) \ar[rrd,bend left] \ar[ddr,bend right] \ar[rd, dashed]\\
	&X_A \ar[r]\ar[d] & X\ar[d,"\pi"]\\
	&\Spec(A)\ar[r] & X_\dR.,
	\arrow["\lrcorner"{anchor=center, pos=0.125}, draw=none, from=2-2, to=3-3]
\end{tikzcd}\]
There exist two presheaves of complete filtered $k$-algebras on $X_\dR$: 
\[\calB_X(A) \defeq \DR(\Spec(A^\red)/X_A),\]
and
\[\bD_{X_\dR}(A)\defeq \DR(\Spec(A^\red)/\Spec(A)).\]
We also naturally get a map $\bD_{X_\dR}\rightarrow \calB_X$. If we forget filtrations, we have
\[(\calB_X)^u \simeq \pi_*\calO_X\]
and
\[(\bD_{X_\dR})^u \simeq \calO_{X_\dR}.\]
It was then established in \cite{CPTVV} that an $n$-shifted Poisson structure on $X$ should be given by a $\bD_{X_\dR}$-linear $n$-shifted Poisson structure on $\calB_X$:
\[\Pois(X,n)\defeq \Pois(\calB_X/\bD_{X_\dR},n).\]
There is a little technical thing to do in order to properly define $\Pois(\calB_X/\bD_{X_\dR},n)$, because one wants the Poisson bracket to be compatible with filtrations. In \cite{CPTVV}, the authors use a trick by passing to ind-objects. In this paper, we establish that $\Pois(\calB_X/\bD_{X_\dR},n)$ should correspond to $\bD_{X_\dR}$-linear $n$-shifted Poisson structures on $\calB_X$ with filtration degree $-2$. Now, by construction, $\calB_X(A)$ and $\bD_{X_\dR}(A)$ are derived foliations over $\Spec(A^\red)$. We prove that, given an $n$-shifted Poisson structure on $X$, we get a compatible family of $A$-linear $(n+1)$-shifted Lagrangian thickenings $(X_A \rightarrow Y_A)$ over $\Spec(A)$ using \cref{prop:relthick}, a relative version of \cite[5.6]{lagthick}. This last proposition relies heavily on the fact that an $n$-shifted Poisson structure on $X$ is a $\bD_{X_\dR}$-linear $n$-shifted Poisson structure on $\calB_X$ of filtration degree $-2$. Explicitly, we get an object in:
\[\lim\limits_{\Spec(A) \rightarrow X_\dR} \Lag\Thick(X_A/\Spec(A),n+1).\]
We prove in \cref{prop:lagthickdesc} that such a family can be glued globally on $X$:
\[\lim\limits_{\Spec(A) \rightarrow X_\dR} \Lag\Thick(X_A/\Spec(A),n+1) \simeq \Lag\Thick(X,n+1).\]
We get that way an $(n+1)$-shifted Lagrangian thickening of $X$.
\subsection*{Other things contained in this paper}
\par Thanks to the tools developed to prove theorem \ref{thm:introaksz}, we also include a few implications:
\begin{thmintro}
	If $X$ is a prestack that has a deformation theory and is locally of finite presentation, then there is an equivalence
	\[\Pois^\nd(X,n)\simeq \Symp(X,n)\]
	between the spaces of non-degenerate $n$-shifted Poisson structures on $X$ and the space of $n$-shifted symplectic structures on $X$.
\end{thmintro}
\par This theorem is already known when $X$ is Artin and locally of finite presentation, but here, we extend it. In particular, it applies to formal neighborhoods of morphisms of Artin stacks. We then give two interesting applications of it. The first one is about moduli of flat connections on a non-necessarily proper smooth scheme $X$. Using \cite{TP2}, we prove that on each $k$-point of
\[\Vect^\nabla(X),\]
the formal neighborhood has a shifted Poisson structure. The second application comes from BV formalism. We sketch a proof that given a function $S : U\rightarrow \bA^1$ on a smooth affine scheme $U=\Spec(A)$, we have a $(-1)$-shifted coisotropic structure on
\[\stCrit(S) \rightarrow \dCrit(S)\]
the inclusion of the strict critical locus into the derived critical locus.
\subsection*{Relation with other work}
\par Another method to do the AKSZ construction is currently under investigation by Melani and Pourcelot (\cite{intpois}). Their strategy relies on integration along fibers of shifted Poisson structures, which is a more general method than what is presented here, we expect that their strategy work on any derived prestack without any extra representability assumptions.
\par We expect that we can recover from our construction known Poisson structures on moduli stacks, for instance, built in \cite{poismod1} or \cite{poismod2}. Comparsion between those two constructions is a non-trivial work and would be interesting to investigate.
\section*{Acknowledgments}
\par I would like to warmly thank my advisor, Damien Calaque, for his constant support and for all our discussions about this work. I also thank Joost Nuiten for an important discussion about formal localization. Finally, I thank Bertrand To\"en for discussions about shifted Poisson structures on moduli stack of flat connections over a non-necessarily proper Deligne--Mumford stack, and for pointing out a gap in an earlier version of this paper.
\section{Preliminaries}
\par In this paper, the word ``category" refers to $\infty$-categories, unless specified otherwise. We fix $k$, a field of zero characteristic, and the following notations:
\begin{itemize}
	\item $\Sp$ is the category of spaces or $\infty$-groupoids,
	\item $\Mod_k$ is the category of $k$-modules, 
	\item $\CAlg_k$ is the category of commutative $k$-algebras, and
	\item if $\calC$ is a stable category with a t-structure, $\calC^\conn$ is the category of \textit{connective} objects of $\calC$.
\end{itemize}
\par We fix $\calC$, any $k$-linear presentably symmetric monoidal category. For $\calO$ an operad in $\calM$ and $\calD$ a category tensored over $\calM$, we denote by
\[\Alg_\calO(\calD)\]
the category of $\calO$-algebras in $\calD$. Recall the usual $k$-linear operads
\begin{itemize}
	\item $\Com$ classifying commutative algebras,
	\item $\Lie$ classifying Lie algebras, and
	\item $\bP_n$ classifying Poisson algebras with a bracket of degree $1-n$.
\end{itemize}
If $\calO = \Com$ is the commutative operad, we will denote by
\[\CAlg(\calD)\]
the category of commutative algebras instead.
\par Recall that a \defterm{derived prestack} is an accessible functor $\CAlg_k^\conn\rightarrow \calS$. We denote by $\dPSt$ the category of derived prestacks. In the rest of this paper, we will refer to such object as ``prestacks" instead, thus omitting the adjective ``derived."
\subsection{Graded, filtered, and complete filtered objects}
\par Let $\Z$ be the poset of integers with the order $\leq$. The category $\calC^{\fil(, \geq 0)}$ of (non-negatively) \defterm{filtered objects} is defined as follows:
\[\calC^{\fil(, \geq 0)} \defeq \Fun\left((\Z_{(\geq 0)},\leq)^{\op},\calC\right).\]
An object $\calC^{\fil}$ will be denoted $F^{\bullet} M$ with maps
\[\dots \rightarrow F^i M \rightarrow F^{i-1}M\rightarrow \dots.\]
An object of $\calC^{\fil, \geq 0}$ will be denoted the same way, but the filtration is bounded on the right:
\[\dots \rightarrow F^1M \rightarrow F^0M.\]
We will denote by $M(i)$ the filtered object given by
\[F^j (M(i)) \defeq \begin{cases} 0 & \text{ if } j < i\text{,}\\ M & \text{elsewhere.}\end{cases}.\]
We set
\[\calC^{\gr(,\geq 0)} \defeq \Fun\left(\Z^{\delta}_{(\geq 0)},\calC\right)\]
for the category of (non-negatively) \defterm{graded objects} of $\calC$, where $\Z^{\delta}$ is the discrete category of integers. An object of $\calC^{\gr(, \geq 0)}$ will be denoted $M^\bullet$. There is a functor denoted $\gr$:
\[\begin{array}{ccc}
	\calC^{\fil(,\geq 0)} & \longrightarrow & \calC^{\gr(, \geq 0)}\\
	F^\bullet M & \longmapsto & \left(\cofib(F^{i+1} M \rightarrow F^i M)\right).
\end{array}\]
that is called the associated graded functor. If $M \in \calC$ and $i \in \Z$, we will denote by $M(i)$ the graded object of $\calC$ with $M$ put in weight $i$ and $0$ elsewhere. The notation is the same as the one in the filtered case; however, it will be clear from the context what object we are talking about. We can see that it is also compatible with the functor $\gr$. $\calC^{\fil(, \geq 0)}$ and $\calC^{\gr(, \geq 0)}$ are $k$-linear and presentable, and $\gr$ commutes with limits and colimits. Using the addition of $\Z$ and $\Z^\delta$, we get presentably symmetric monoidal products on both of them by Day convolution (\cite[2.2.6.]{HA}).
\par Recall that an object $F^{\bullet} M \in \calC^{\fil}$ is called \defterm{complete} if $\lim\limits_n F^n M \simeq 0$. We denote $\calC^{\cfil}$ and $\calC^{\cfil,\geq 0}$ the full subcategories spanned by complete filtered objects. Remark that $M(i)$ is always complete by definition. There is an adjunction:
\[\widehat{(\vari)} : \calC^{\fil(, \geq 0)} \rightleftarrows \calC^{\cfil(, \geq 0)} : i\]
where $i$ is the inclusion and $\widehat{(\vari)}$ is the completion. It is actually a left exact localization at morphisms inducing equivalences on associated graded. $\calC^{\cfil(, \geq 0)}$ is then $k$-linear presentably symmetric monoidal (the symmetric monoidal product being given by the completion of the monoidal product of filtered objects). The functor $\gr$ factors as
\[\gr : \calC^{\cfil(, \geq 0)} \rightarrow \calC^{\gr(,\geq 0)}.\]
It is conservative, symmetric monoidal, and has right and left adjoints. It induces functors
\[\CAlg(\calC^{(\wedge)\fil(,\geq 0)}) \rightarrow \CAlg(\calC^{\gr(, \geq 0)})\]
that also have right and left adjoints (in the complete case, the induced functor is conservative). It is also well-behaved with the internal Hom functors in the complete case (see \cite[2.1.]{lagthick}):
\begin{lm}\label{lm:grhom}
	Let $A \in \CAlg(\calC^{\cfil})$ and $X,Y \in \Mod_A(\calC^{\cfil})$. There is a natural equivalence
	\[\gr\underline\Hom_A(X,Y) \simeq \underline\Hom_{\gr A}(\gr X,\gr Y).\]
\end{lm}
\par If $G : \calC \rightarrow \calD$ is lax monoidal, where $\calD$ satisfies the same assumptions as $\calC$, then we have induced lax monoidal functors
\[G^{\fil(, \geq 0)} : \calC^{\fil(, \geq 0)} \rightarrow \calD^{\fil(, \geq 0)}\]
and
\[G^{\gr(, \geq 0)} : \calC^{\gr( \geq 0)} \rightarrow \calD^{\gr(, \geq 0)}.\]
If it is, moreover, symmetric monoidal and commutes with colimits, then the induced lax structure is symmetric monoidal because the Day convolution tensor product is a left Kan extension. If $G$ commutes with limits, there is an induced lax monoidal functor
\[G^{\cfil(, \geq 0)} : \calC^{\cfil(, \geq 0)} \rightarrow \calD^{\cfil(, \geq 0)}.\]
There is an ``underlying object" symmetric monoidal colimit-preserving functor
\[(\vari)^u : \calC^{\fil(, \geq 0)} \rightarrow \calC\]
sending $F^\bullet M$ to $\colim\limits_n F^n M$. It also induces a left adjoint on the level of algebras:
\[(\vari)^u : \CAlg(\calC^{\fil(, \geq 0)}) \rightarrow \CAlg(\calC).\]
Note that in the non-negatively filtered case, these functors send $F^\bullet M$ to $F^0 M$. In the complete case, it is slightly different. The underlying object functor is neither a right adjoint nor a left adjoint. However, the functor $F^0(\vari)$ is lax monoidal and has a symmetric monoidal left adjoint given by $(\vari)(0)$. By tradition (where complete filtered objects are seen as mixed graded objects instead in shifted symplectic geometry), we will denote $F^0(\vari)$ by $|\vari|$. When we apply $|\vari|$ to an object, we say that we ``realize" it. The functor $(\vari)(0)$ also has a left adjoint traditionally called ``left realization." Explicitly, it sends $F^\bullet M$ to $\cofib\left(F^1 M \rightarrow (F^\bullet M)^u\right)$.
\par When an object is doubly filtered or graded, we will need to distinguish notations. In that case, the notation $\langle n \rangle$ will play the same role as $(n)$, but one will be applied to the first filtration/grading and the other one to the second filtration/grading. Usually we choose $\langle n \rangle$ for filtrations coming from Hodge filtrations of de Rham forms, and $(n)$ will be used for internal filtration, where we consider our base objects to be naturally given with an extra filtration. 
\subsection{de Rham forms}
\par We refer to \cite[Sections 2 and 3.1.]{lagthick} for definitions and basic properties about the de Rham functor. Recall that there is a functor
\[\DR(\vari/\vari) : \CAlg(\calC)^{\Delta^1} \rightarrow \CAlg^{\cfil}\]
that sends $B \rightarrow A$ to $\DR(A/B)$, the complete complex of $B$-linear de Rham forms on $A$ with the Hodge filtration. There is a version for prestacks $X \rightarrow Y$, where
\[\DR(X/Y) \in \QCoh(Y)^{\cfil}.\]
\subsection{Shifted symplectic and Lagrangian structures}
\par We refer to \cite[3.1.]{lagthick} for a review on the theory of shifted symplectic and Lagrangian structures for prestacks. Recall that for $A\in \CAlg(\calC)$ with a perfect cotangent complex, and $\calL \in \calC$ an invertible object, we denote by
\[\Symp(A,\calL)\]
the subspace of $\rm{pre}\Symp(A,\calL) \defeq \Map_{\Mod(\calC^{\cfil})}(\mathbf{1}[-2]\langle 2 \rangle,\DR(A)\otimes\calL)$ consisting of $\calL$-twisted symplectic structures on $A$, and if $f : B \rightarrow A$ is a morphism in $\CAlg(\calC)$ with a perfect cotangent complex, we denote by
\[\Lag(f,\calL)\]
the subspace of the pullback
\[\begin{tikzcd}
	\Isot(f,\calL) \ar[r]\ar[d]& \Map_{\Mod(\calC^{\cfil})}(\mathbf{1}[-2]\langle 2 \rangle,\DR(B)\otimes\calL)\ar[d] \\
	* \ar[r,"0"] & \Map_{\Mod(\calC^{\cfil})}(\mathbf{1}[-2]\langle 2 \rangle,\DR(A)\otimes\calL)
	\arrow["\lrcorner"{anchor=center, pos=0.125}, draw=none, from=1-1, to=2-2]
\end{tikzcd}\]
consisting of $\calL$-twisted Lagrangian structures on $f$ (it induces a symplectic structure on $B$). If $g : C \rightarrow B$ is another morphism, we can see $f$ as a morphism in $\Mod_C(\calC)$. Then we can take $\calL \in \Mod_C(\calC)$ invertible, and we get relative versions
\begin{center}\begin{tabular}{cc}
	$\rm{pre}\Symp(B/C,\calL)$, & $\Symp(B/C,\calL)$,\\\\
	$\Isot(f/C,\calL)$, and & $\Lag(f/C,\calL)$.
\end{tabular}\end{center}
Now, recall from \cite[Section 2]{sympgrpd} that if $f : X \rightarrow Y$ is a morphism of prestacks having cotangent complexes, we can define the spaces
\begin{center}\begin{tabular}{cc}
	$\rm{pre}\Symp(X,\calL)$, & $\Symp(X,\calL)$,\\\\
	$\Isot(f,\calL)$, and & $\Lag(f,\calL)$
\end{tabular}\end{center}
of respectively $n$-shifted (pre)symplectic structures on $X$ and $n$-shifted isotropic and Lagrangian structures on $f$. If $g : Y \rightarrow Z$ is another morphism of prestacks, and $X$ and $Y$ are assumed to only have a cotangent complex relative to $Z$, there are also relative versions
\begin{center}\begin{tabular}{cc}
	$\rm{pre}\Symp(X/Z,\calL)$, & $\Symp(X/Z,\calL)$,\\\\
	$\Isot(f/Z,\calL)$, and & $\Lag(f/Z,\calL)$.
\end{tabular}\end{center}
\subsection{Formal geometry}
\par Recall that a morphism $f : B \rightarrow A$ of commutative $k$-algebra $A$ is said to be
\begin{itemize}
	\item \defterm{almost of finite presentation} if $A$ is almost compact as a $B$-algebra, and
	\item \defterm{a nilpotent extension} if $\pi_0(B)\rightarrow\pi_0(A)$ is surjective with nilpotent kernel.
\end{itemize}
The ring $A$ is said to be \defterm{almost of finite presentation} if $k \rightarrow A$ is almost of finite presentation.
 A morphism of prestacks $X\rightarrow Y$ is said to be \defterm{convergent} if for any $A\in \CAlg_k^\conn$, the following diagram:
\[\begin{tikzcd}
	X(A) \ar[r] \ar[d] & \lim X(\tau_{\leq n}A)\ar[d]\\
	Y(A) \ar[r] & \lim Y(\tau_{\leq n} A)
\end{tikzcd}\]
is Cartesian. The prestack $X$ is said to be \defterm{convergent} if $X \rightarrow *$ is convergent. We denote by $\dPSt^\conv$ the full subcategory of $\dPSt$ composed of convergent prestacks, and by $X^\conv \defeq \lim X(\tau_{\leq n}(\vari))$ the \defterm{convergent completion} of $X$. It is said to have a \defterm{deformation theory} if it is convergent and for any pullback square
\[\begin{tikzcd}
	B'\ar[r] \ar[d] & B\ar[d]\\
	A'\ar[r] & A
	\arrow["\lrcorner"{anchor=center, pos=0.125}, draw=none, from=1-1, to=2-2]
\end{tikzcd}\]
of connective algebras with $A'\rightarrow A$ nilpotent,
\[X(B') \simeq X(B)\underset{X(A)}\times X(A').\]
\par A map $X\rightarrow Y$ of derived prestack is said to be \defterm{locally almost of finite presentation} of \defterm{lafp} for short if it is convergent and for any diagram $I \rightarrow \CAlg_k^{\conn}$ of $n$-truncated connective $k$-algebras, whose colimit is $A$, the following square:
\[\begin{tikzcd}
	\colim X(A_i) \ar[r] \ar[d] & X(A)\ar[d]\\
	\colim Y(A_i) \ar[r] & Y(A)
\end{tikzcd}\]
is Cartesian. The prestack $X$ is said to be lafp if $X\rightarrow *$ is lafp. We denote by $\dPSt^\lafp$ the full subcategory of $\dPSt$ composed of lafp prestacks.
\begin{warning}
	Being locally of finite presentation implies being locally almost of finite presentation, but the converse is false. For instance, the free algebra $k[x_i]$ for $x_i$ of degree $-2i$ is (locally) almost of finite presentation but not (locally) of finite presentation.
\end{warning}
We also define the following prestacks:
	\[X_\dR(A) \defeq X(A^\red),\]
called the \defterm{de Rham prestack} of $X$, and
	\[X_\red \defeq \colim\limits_{\Spec(A)\rightarrow X} \Spec(A^\red).\]
called the \defterm{reduced prestack} associated to $X$, note that $(\vari)_\red$ is left adjoint to $(\vari)_\dR$. For $f : X \rightarrow Y$ a morphism of prestacks, we set:
	\[Y^\wedge_X \defeq X_\dR \underset{Y_\dR}\times Y.\]
to be \defterm{the formal completion} of $X$ along $f$.
\begin{lm}\label{lm:drdefth}
	Let $X$ be any prestack, then $X_\dR$ has a deformation theory.
\end{lm}
\begin{proof}
	If $A \rightarrow A'$ is a nilpotent extension, then $A^\red \simeq A'^\red$, so $X_\dR$ sends Cartesian diagrams:
	\[\begin{tikzcd}
		B' \ar[r] \ar[d] & B \ar[d]\\
		A'\ar[r] & A
		\arrow["\lrcorner"{anchor=center, pos=0.125}, draw=none, from=1-1, to=2-2]
	\end{tikzcd}\]
	into Cartesian diagrams. It is also convergent because
	\[\lim\limits_n X_\dR(\tau_{\leq n} A) \simeq \lim\limits_n X((\tau_{\leq n} A)^\red) \simeq \lim\limits_n X(A^\red) \simeq X_\dR(A).\]
\end{proof}
Moreover, mapping prestacks usually preserves the fact of having a deformation theory:
\begin{prop}\label{prop:mapdef}
	Let $Y$ be a prestack and $X$ a prestack with deformation theory. Then
	\[\MapSt(Y,X)\]
	has a deformation theory.
\end{prop}
\begin{proof}
	We can reduce to the affine case since $\MapSt(Y,X) \simeq \lim\limits_{\Spec(A) \rightarrow Y} \MapSt(\Spec(A),X)$ and having a deformation theory is stable by taking limits. In that case, we can follow Lurie's proof of \cite[19.1.3.1. (2)-(3)]{SAG}. He uses the fact that, $B$-points of $\MapSt(\Spec(A),X)$ are $B\otimes A$-points of $X$, and $B \otimes \vari$ preserves eventually constant sequences on each $\pi_n$, pullbacks of $k$-algebras along surjective morphisms, and nilpotent extensions.
\end{proof}
\par In general, being lafp is not preserved by taking mapping prestacks. However, it is possible to enforce this condition and prove that it preserves all the good properties we want on our prestack. Because being lafp is stable by colimit of convergent prestacks, we have a right adjoint of the inclusion
\[\dPSt^\lafp \rightarrow \dPSt^\conv\]
between convergent and lafp morphisms. We denote this right adjoint $(\vari)^\lafp$.
\begin{prop}\label{prop:lafp} The following are satisfied:
	\begin{enumerate}
		\item If $X$ is a prestack that has a deformation theory, then $X^\lafp$ has a deformation theory.
		\item If $f : Y \rightarrow X$ is a morphism of convergent prestacks that has a perfect cotangent complex, and if $f^\lafp$ has a perfect cotangent complex, then for any $x : \Spec(B) \rightarrow Y^\lafp$ with $B$ almost of finite presentation,
		\[x^*\iota_Y^* \cot_f \simeq x^*\cot_{f^\lafp}\]
		where $\iota_Y : Y^\lafp \rightarrow Y$ is the counit of the adjunction.
	\end{enumerate}
\end{prop}
\begin{proof}Let $X$, and $f : Y \rightarrow X$ be as in the proposition.
	\begin{enumerate}
		\item Asume $X$ has a deformation theory, because $X^\lafp$ is lafp, to check that it has a deformation theory, it suffices to check that for any pullback diagrams
		\[\begin{tikzcd}
			B'\ar[d]\ar[r] & B\ar[d]\\
			A'\ar[r] & A
			\arrow["\lrcorner"{anchor=center, pos=0.125}, draw=none, from=1-1, to=2-2]
		\end{tikzcd}\]
		of almost of finite presentation $k$-algebras with $A' \rightarrow A$ a nilpotent extension, there is a unique lift:
		\[\begin{tikzcd}
			\Spec(A) & \Spec(A') \\
			\Spec(B) & \Spec(B') \\
			&& X^\lafp
			\arrow[from=1-1, to=1-2]
			\arrow[from=1-1, to=2-1]
			\arrow[from=1-2, to=2-2]
			\arrow[bend left, from=1-2, to=3-3]
			\arrow[from=2-1, to=2-2]
			\arrow[bend right, from=2-1, to=3-3]
			\arrow[dashed, from=2-2, to=3-3]
		\end{tikzcd}\]	
		but because $(\vari)^\lafp$ is a right adjoint, this amounts to checking that there is a unique lift:
		\[\begin{tikzcd}
			\Spec(A) & \Spec(A') \\
			\Spec(B) & \Spec(B') \\
			&& X
			\arrow[from=1-1, to=1-2]
			\arrow[from=1-1, to=2-1]
			\arrow[from=1-2, to=2-2]
			\arrow[bend left, from=1-2, to=3-3]
			\arrow[from=2-1, to=2-2]
			\arrow[bend right, from=2-1, to=3-3]
			\arrow[dashed, from=2-2, to=3-3]
		\end{tikzcd}\]
		This is true since $X$ has a deformation theory.
	\item Because $B$ is almost of finite presentation, we know that $X(B) \simeq X^\lafp(B)$. Let $M$ be a perfect $B$-module, we know that $B[M]$ is still almost of finite presentation and:
		\begin{align*}
			\Map(x^*\iota_Y^* \cot_f,M) &\simeq \Map_{\Spec(B)/}(\Spec(B[M]),Y)\underset{\Map(\Spec(B[M]),X)}\times \{0\}\\
						    &\simeq \Map_{\Spec(B)/}(\Spec(B[M]),Y^\lafp)\underset{\Map(\Spec(B[M]),X^\lafp)}\times \{0\}\\
						    &\simeq \Map(x^*\cot_{f^\lafp},M)
		\end{align*}
		because $\cot_f$ and $\cot_{f^\lafp}$ are perfect, we conclude that $x^*\iota_Y^*\cot_f \simeq x^*\cot_{f^\lafp}$.
	\end{enumerate}
\end{proof}
\begin{defn}
	A prestack $X$ is said to be \defterm{formal} if it has a deformation theory, has a perfect cotangent complex, and is locally of finite presentation.
\end{defn}
\begin{remark}
	In \cite{CPTVV}, the authors call a formal prestack a prestack that has a deformation theory and is lafp. Here we put a more restrictive assumption, asking for the prestack to be locally of finite presentation and have a perfect cotangent complex instead. It is because it is the class of prestacks where shifted Poisson geometry behaves and interracts well with shifted symplectic geometry. 
\end{remark}
\begin{warning}
	It is not true that a prestack $X$ having a deformation theory and being locally of finite presentation has automatically a perfect cotangent complex. For instance, take $X$ to be the prestack left Kan extended from
	\[A \mapsto \Mod_A^\conn\]
	for $A$ of finite presentation. It is locally of finite presentation and has a deformation theory thanks to \cite[16.2.0.2.]{SAG}. Take an $A$-point $M$, for $A$ of finite presentation. Then by \cite[16.5.4.1.]{SAG}, the relative tangent complex at $M$ is $\mathrm{End}(M)$ which is neither perfect nor connective.
\end{warning}

\begin{defn}
	A map $X\rightarrow Y$ of prestacks is a \defterm{formal thickening} if
	\begin{itemize}
		\item $Y$ has a deformation theory,
		\item $X \rightarrow Y$ is lafp,
		\item and the induced map $X_\dR \rightarrow Y_\dR$ is an equivalence.
	\end{itemize}
	We let $\Thick(X)$ be the category of \defterm{formal thickenings} of $X$.
\end{defn}
\par Later, we will consider Lagrangian and isotropic thickenings. We introduce the notations
\[\Lag\Thick(X,n)\]
and
\[\Isot\Thick(X,n)\]
for the spaces of thickenings of $X$ with an $n$-shifted Lagrangian structure (resp. an $n$-shifted isotropic structure).
\begin{defn}
	A \defterm{derived foliation} over $A\in \CAlg_k$ is the data of $\calF \in \CAlg^{\cfil}_k$, an equivalence $\gr^0 \calF\simeq A$, and such that the natural morphism
	\[\Sym_{\gr^0 \calF}(\gr^1\calF\langle 1 \rangle) \rightarrow \gr\calF\]
	is an equivalence. We set $\cot_\calF \defeq \gr^0\calF[1]$ and call it the \defterm{cotangent complex} of $\calF$. A derived foliation is called \defterm{perfect} (resp. \defterm{almost perfect} if $\cot_\calF$ is perfect (resp. almost perfect) as an $A$-module; in that case, we set $\tan_\calF\defeq \cot_\calF^\vee$. We denote by $\Fol_A$ the category of derived foliations over $A$, and $\Fol^\ap_A$ the full subcategory of $\Fol_A$ consisting of almost perfect derived foliations. To avoid ambiguity, we denote by $\cot^\int_\calF$ the cotangent complex of $\calF$ seen as a complete filtered commutative algebra, in that case, it will be a complete filtered $\calF$-module.
\end{defn}
Recall from \cite[3.36.]{lagthick} that, when $A$ is almost of finite presentation, we can deduce from \cite{follie} and \cite{deffol} an integration functor
\[\Fol_A^{\ap,\op} \rightarrow \Thick(\Spec(A)).\]
For $\calF \in \Fol_A^\ap$, we denote by $[\Spec(A)/\calF]$ the image of this functor in $\Thick(\Spec(A))$. It is called the \defterm{prestack of leaves} of $\calF$. By \cite[12.5.3.]{SAG} and \cite[5.3.]{deffol}, $[\Spec(A)/\calF]$ will always be of the form $(\colim \Spec(A_i))^\conv$, where $A_i \rightarrow A$ are Artin extensions, and the colimit is sifted.
\begin{warning}
	Our definition of the category of derived foliations is rather unconventionnal: in \cite{livrefol}, the authors define the category of derived foliations over $A$
	\[\mathcal{F}ol(A/k)\]
	to be the opposite category $\Fol_A^\op$. We have decided to take a more algebraic definition, we see a derived foliation $\calF$ directly as its de Rham complex denoted $\DR(\calF)$ in \cite{livrefol}, and not as its formal spectrum (formal in the categorical sense).
\end{warning}
We can freely talk about shifted symplectic structures on this prestack:
\begin{prop}[{\cite[3.39.]{lagthick}}]
	If $\calF$ is perfect and $A$ is of finite presentation, then
	\[[\Spec(A)/\calF]\]
	has a perfect cotangent complex, and $\cot^\int_\calF$ is perfect.
\end{prop}
\par We can compute internal relative cotangent complexes:
\begin{lm}\label{lm:relcot}
	Let $\calF \rightarrow \calG$ be a morphism of derived foliations over $A$, then
	\[\gr\cot^\int_{\calG/\calF} \simeq \cofib(\cot_\calF \rightarrow \cot_\calG)\underset{A}\otimes\gr\calG[-1](1).\]
\end{lm}
\begin{proof}
	We have a fiber sequence
	\[\cot_{\gr\calF}\underset{\gr\calF}\otimes\gr\calG \rightarrow \cot_{\gr\calG} \rightarrow \cot_{\gr\calG/\gr\calF}.\]
	that gives a fiber sequence
	\[\left(\cot_A(0) \oplus \cot_\calF[-1](1)\right)\underset{A}\otimes\gr\calG \rightarrow \left(\cot_A(0) \oplus \cot_\calG[-1](1)\right)\underset{A}\otimes \gr\calG \rightarrow\cot_{\gr\calG/\gr\calF}\]
	so that
	\[\gr\cot^\int_{\calG/\calF}\simeq \cot_{\gr\calG/\gr\calF}\simeq \cofib(\cot_\calF \rightarrow \cot_\calG)\underset{A}\otimes\gr\calG[-1](1).\]
\end{proof}
\section{Formal localization}\label{sec:formloc}
\par Fix $X$ a formal prestack. Recall from \cite[2]{CPTVV} the technique called formal localization. Roughly speaking, $X$ is ``almost affine" over its de Rham prestack. If
\[\begin{tikzcd}
	X_A \ar[r] \ar[d] & X \ar[d]\\
	\Spec(A) \ar[r] & X_\dR
	\arrow["\lrcorner"{anchor=center, pos=0.125}, draw=none, from=1-1, to=2-2]
\end{tikzcd}\]
is a pullback diagram, we can check by abstract nonsense \cite[2.1.8.]{CPTVV} that the prestack $X_A$ is a formal thickening of $\Spec(A^\red)$ when $A$ is almost of finite presentation. Because $X$ is locally of finite presentation, $X_\dR$ is locally of finite presentation as well, and we can take $A$ to be of finite presentation. In that case, $X_A$ is locally of finite presentation and will be of the form $[\Spec(A^\red)/\calF_A]$, where $\calF_A \defeq \DR(\Spec(A^\red)/X_A)$ is the derived foliation associated to $\Spec(A^\red) \rightarrow X_A$. We have a family over the de Rham prestack
\[\calB_X : \begin{array}{rcl}\dAff^\op_{X_\dR} & \longrightarrow &\CAlg_k^{\cfil}\\
	\Spec(A) & \longmapsto & \DR(\Spec(A^\red)/X_A)\end{array}\]
called the \defterm{filtered principal parts} of $X$. We can also see $\Spec(A)$ as a formal thickening of $\Spec(A^\red)$, we then set
		\[\bD_{X_\dR} : \begin{array}{rcl} \dAff^\op_{X_\dR} & \longrightarrow & \CAlg_k^{\cfil}\\
	\Spec(A) & \longmapsto & \DR(A^\red/A)\end{array}\]
to be the \defterm{filtered crystalline structure presheaf} of $X$. There is a functor
\[\bD_{X_\dR} \rightarrow \calB_X\]
giving to $\calB_X$ a $\bD_{X_\dR}$-linear structure. It is possible to compute forms on $X$ in terms of forms on $\calB_X$. To state and prove the precise statement, we need to generalize \cite[3.43.]{lagthick} for relative forms.
\subsection{Forms on derived foliations}
We fix $B$ a commutative $k$-algebra of finite presentation, and $X\defeq \Spec(B)$. This is a generalization of \cite[3.43.]{lagthick}:
\begin{thm}\label{thm:formformal}
	Let $\calF \rightarrow \calG$ be a morphism of almost perfect derived foliations over $B$. Then we have:
	\[\Gamma\left([X/\calF],\DR([X/\calG]/[X/\calF])\right) \simeq |\DR^\int(\calG/\calF)|^\int.\]
\par 	Moreover, we have a commuting diagram:
	\[\begin{tikzcd}
		\left|\Sym_\calG(\cot^\int_{\calG/\calF}[-1]\langle 1 \rangle)\right|^\int \ar[d] & \ar[l,"\sim"'] \Gamma\left([X/\calG],\Sym_{[X/\calG]}(\cot_{[X/\calG]/[X/\calF]}[-1]\langle 1 \rangle)\right)\ar[d]\\
		\left|\Sym_{B(0)}(\cot^\int_{\calG/\calF}\underset\calG\otimes B(0)[-1]\langle 1 \rangle)\right|^\int & \ar[l,"\sim"'] \Gamma\left(X,\Sym_{\calO_X}(\pi^*\cot_{[X/\calG]/[X/\calF]}[-1]\langle 1 \rangle)\right)
	\end{tikzcd}\]
\end{thm}
\par To prove this theorem, we follow a strategy similar to the one used in the proof of \cite[3.43.]{lagthick}. We proceed by proving several lemmas.
\begin{lm}\label{lm:lm1}
	Let $Y \rightarrow Z$ be a map of prestacks, then
	\[\Gamma(Z,\DR(Y/Z)) \simeq \Gamma(Y,\DR(Y^\conv/Z^\conv))\]
	as complete filtered commutative algebras.
\end{lm}
\begin{proof}
	For now, we will omit all occurences of global sections, remembering that our object always live in $\Mod_k$ instead of $\QCoh(Z)$ (or $\QCoh(Z^\conv)$).\\
	By definition,
	\[\DR(Y/Z)\simeq \lim_{\Spec(A) \rightarrow Z} \DR(Y_A/\Spec(A)).\]
	And,
	\[\DR(Y/Z) \simeq \lim_{\Spec(A) \rightarrow Z}\lim_{\Spec(C) \rightarrow Y_A} \DR(C/A).\]
	We claim that we can restrict the diagram where this limit is taken to the subdiagram consisting of those $A$ and $C$ that are eventually coconnective. First, $\DR(\vari/A)$ is convergent, the proof goes the same way as the proof of \cite[3.1.]{cohohall3CY}, so we can restrict to $C$ that are eventually coconnective. Now, we prove that, when $C$ is eventually coconnective,
	\[\DR(C/A) \simeq \lim_{n\gg 0} \DR(C/A_n)\]
	where $A_n \defeq \tau^{\geq -n} A$ is the $n$-truncation of $A$, $n$ is picked to be big enough so that a factorization $A \rightarrow A_n \rightarrow C$ exists. By definition, there is a natural map
	\[\DR(C/A) \rightarrow \lim_{n\gg 0} \DR(C/A_n).\]
	On the associated graded, this map is given by
	\[\Sym^i_C(\cot_{C/A}[-1]) \rightarrow \lim \Sym^i_C(\cot_{C/A_n}[-1]).\]
	Thanks to the exact sequence
	\[\cot_{A_n/A}\otimes_{A_n} C \rightarrow \cot_{C/A} \rightarrow \cot_{C/A_n},\]
	the fact that $\cot_{A_n/A}$ is $(n+2)$-connected (thanks to \cite[7.4.3.2]{HA} and $A \rightarrow A_n$ has an $(n+2)$-connective cofiber), and the fact that pullback is right $t$-exact, we have that
	\[\cot_{C/A}\rightarrow \cot_{C/A_n}\]
	is $(n+1)$-connective. The same argument yields to the fact that the tower
	\[\dots \rightarrow \cot_{C/A_{n+1}} \rightarrow \cot_{C/A_n} \rightarrow \dots\]
	is eventually constant on each homotopy groups. By the Milnor exact sequence, for any $n$, $\pi_m \lim_n \cot_{C/A^n} \simeq \lim_n \pi_m \cot_{C/A^n}$, and $\cot_{C/A} \simeq \lim \cot_{C/A^n}$. The same argument works for higher $\Sym$s because the fiber of
	\[\Sym^i_C(\cot_{C/A}[-1]) \rightarrow \Sym^i_C(\cot_{C/A^n}[-1])\]
	will be $(n+2-i)$-connective (like in the proof of \cite[3.1.]{cohohall3CY}). So, we may rephrase our description of $\DR(Y/Z)$ as follows:
	\[\DR(Y/Z) \simeq \lim_{\Spec(A) \rightarrow Z}\lim_{\underset{C \text{ eventually coconnective}}{\Spec(C) \rightarrow Y_A}}\lim_{\underset{A' \text{ eventually coconnective}}{A \rightarrow A' \rightarrow C}} \DR(C/A').\]
	So
	\begin{align*}\DR(Y/Z) &\simeq \lim_{\underset{A' \text{ eventually coconnective}}{\Spec(A) \rightarrow Z}}\lim_{\underset{C \text{ eventually coconnective}}{\Spec(C) \rightarrow Y_A}} \DR(C/A)\\
		&\simeq \DR(Y^\conv/Z^\conv).\end{align*}
	
\end{proof}
\par This is the ``affine" case of the theorem:
\begin{lm}\label{lm:lm2}
	Let $A \rightarrow B$ and $C \rightarrow B$ two Artin extensions, with a map $A \rightarrow C$, then
	\[\DR(C/A) \simeq |\DR^\int(\DR(B/C)/\DR(B/A))|^\int.\]
\end{lm}
\begin{proof}
	To construct the natural map $\DR(C/A) \rightarrow |\DR^\int(\DR(B/C)/\DR(B/A))|^\int$, one has to compute the underlying commutative algebra of $|\DR^\int(\DR(B/C)/\DR(B/A))|^\int$. Recall from \cite[3.45.]{lagthick}, that $\DR(B/A)$ and $\DR(B/C)$ are respectively equivalent to $\adic(A\rightarrow B)$ and $\adic(C\rightarrow B)$. The functor $\adic(\vari)$ preserves being almost of finite presentation, so 
	\[\adic(A\rightarrow B) \rightarrow \adic(C\rightarrow B)\]
	is almost of finite presentation. By \cite[2.18. (1)]{deffol}, the underlying morphism of commutative algebra is $A \rightarrow C$, an animated surjection. We can apply again an analog of \cite[2.18. (1)]{deffol} to the filtered world (in that case, propositions that are used to prove it need to be adapted to the case where the ideal $I$ is filtered), to get that $\adic^\int(\adic(A\rightarrow B) \rightarrow \adic(C\rightarrow B))$ is complete. We have:
	\begin{align*}
	||\DR^\int(\adic(C\rightarrow B)/\adic(A\rightarrow B))|^\int |^\ext &\simeq ||\DR^\int(\adic(C\rightarrow B)/\adic(A\rightarrow B))|^\ext |^\int\\
									     &\simeq ||\adic^\int(\adic(C\rightarrow B)/\adic(A\rightarrow B))|^\ext |^\int
\\
									     &\simeq|\adic(A\rightarrow B)| \\
									     &\simeq A.
	\end{align*}
	When we compute the $0$th graded piece of this we end up on $C$ for a similar reason. This gives a map
	\[\DR(C/A) \rightarrow |\DR^\int(\DR(B/C)/\DR(B/A))|^\int.\]
	Unwinding classical computations on the associated graded, we obtain the morphism:
	\[\Sym(\cot_{C/A}[-1]) \rightarrow |\Sym(\cot^\int_{\DR(B/C)/\DR(B/A)}[-1])|^\int\]
	This is an equivalence by a similar argument as in the proof of \cite[3.46.]{lagthick}.
\end{proof}
The last lemma states that the internal de Rham complex behaves well with some cosifted limits.
\begin{lm}\label{lm:lm3}
	If $\calF\rightarrow\calG$ is a morphism of almost perfect derived foliations, such that $\calF \simeq \lim \calF_i$ is a cosifted limit, $\calF_i \in \Fol^\ap_B$, such that $\cot_\calF \simeq \lim \cot_{\calF_i}$ and $\cot_\calG \simeq \lim \cot_{\calG\otimes_\calF \calF_i}$, and $\calG \otimes_\calF \calF_i \simeq \lim \calG_{ij}$ is a cosifted limit with the same assumption on the cotangent complex, then
	\[\DR^\int(\calG/\calF)\simeq \lim_i \lim_j \DR^\int(\calG_{ij}/\calF_i).\]
\end{lm}
\begin{proof}
	If we take the double associated graded of the natural map
	\[\DR^\int(\calG/\calF) \rightarrow \lim_i \lim_j \DR^\int(\calG_{ij}/\calF_i)\]
	and unwind it, we get on the left hand side
	\[\Sym_B(\cofib(\cot_\calF \rightarrow \cot_\calG)\langle 1 \rangle (1) [-2] \oplus \cot_\calG(1)[-1])\]
	and on the right hand side
	\[\lim_i\lim_j \Sym_B(\cofib(\cot_{\calF_i} \rightarrow \cot_{\calG_{ij}})\langle 1 \rangle (1) [-2] \oplus \cot_{\calG_{ij}}(1)[-1])\]
	We can then put limits in the $\Sym$ thanks to \cite[3.34.]{follie}.
\end{proof}
\begin{proof}[Proof of \cref{thm:formformal}.]
	If $Y = [X/\calG]$ and $Z = [X/\calF]$, we have a map $Y \rightarrow Z$. We know that, by definition, 
	\[Z \simeq (\colim_i \Spec(B_i))^\conv\]
	where $B_i \rightarrow B$ are Artin extensions. In that case, we put $Y_i \defeq \Spec(B_i) \times_Z Y$. The prestack $Y_i$ will be of the form $[X/(\calG\otimes_\calF \DR(B/B_i)]$ ($Y_i$ is a formal thickening of $X$ so we can compare it with $[X/(\calG\otimes_\calF \DR(B/B_i)]$ by comparing their cotangent complexes, which turns out to be the same). So we can pick $C_{ij}$ Artin extensions of $B$ such that $Y_i \simeq (\colim_j \Spec(C_{ij}))^\conv$, so
	\[Y \simeq (\colim_{ij} \Spec(C_{ij}))^\conv\]
	by universality of colimits. On the level of foliations, we have that
	\[\calF \simeq \lim \DR(B/B_i)\]
	with $\cot_\calF \simeq \lim \cot_{B/B_i}$ (cf. \cite[3.18.]{lagthick}, the same goes for $\calG\otimes_\calF \DR(B/B_i)$ and $\calG$. Now, we have
	\begin{align*}
		\Gamma(Z,\DR(Y/Z)) &\simeq \Gamma(Z,\DR((\colim_{ij}\Spec(C_{ij}))^\conv/(\colim_i \Spec(B_i))^\conv))\\
			 &\simeq \Gamma(\colim_i \Spec(B_i),\DR(\colim_{ij}\Spec(C_{ij})/\colim_i \Spec(B_i))) & \text{by \cref{lm:lm1},}\\
			 &\simeq \lim_i \lim_j \DR(C_{ij}/B_i)\\
			 &\simeq \lim_i \lim_j |\DR^\int(\DR(B/C_{ij})/\DR(B/B_i))|^\int & \text{by \cref{lm:lm2},}\\
			 &\simeq |\DR^\int(\calG/\calF)|^\int & \text{by \cref{lm:lm3}.}
	\end{align*}
	The proof for the commutative diagram is done the same way as in \cite[3.43]{lagthick}.
\end{proof}
\begin{remark}
	It is legit to ask if \cref{thm:formformal} remains true before taking global sections. We believe it is true but we have to make sense of the $[X/\calF]$-linear structure of the right hand side. In fact, this equivalence should make sense in $\rm{ProCoh}([X/\calF])$.
\end{remark}
A particular case of the latter is the following:
\begin{thm}\label{thm:descdR}
	There is an equivalence of complete filtered algebras
	\[\DR(X) \simeq \left|\Gamma(X_\dR,\DR(\calB_X/\bD_{X_\dR}))\right|^\int.\]
\end{thm}
\begin{proof}
	This was already proven when $X$ is Artin locally of finite presentation in \cite[2.4.12.]{CPTVV}. We give a proof when $X$ is only assumed to be formal. For the first equivalence, we have
	\[\DR(X) \simeq \DR(X/X_\dR) \simeq \lim\limits_{\Spec(A) \rightarrow X_\dR} \DR(X_A/\Spec(A))\]
	Because $X$ is locally of finite presentation, $X_\dR$ is locally of finite presentation as well, and so we can restrict this limit to $k$-algebras $A$ of finite presentation. It suffices to check that $\DR(X_A/\Spec(A)) \simeq \left|\DR^\int(\calB_X(A)/\bD_{X_\dR}(A))\right|^{\int}$. By \cref{thm:formformal} applied to $\bD_{X_\dR}(A)\rightarrow\calB_X(A)$ and because $[\Spec(A^\red)/\calB_X(A)]\simeq X_A$ and $[\Spec(A^\red)/\bD_{X_\dR}(A)]\simeq \Spec(A)$ we get our equivalence.
\end{proof}
\par Thanks to this theorem, we can relate shifted Lagrangian structures to internal shifted Lagrangian structures:
\subsection{Shifted Lagrangian structures on formal thickenings}
\par Let $\calF \rightarrow \calG\overset{f}\rightarrow \calH$ be two morphisms of almost perfect derived foliations over $B$ of finite presentation, assume moreover that the fibers of $\cot_\calF\rightarrow\cot_\calG$ and $\cot_\calF \rightarrow \cot_\calH$ are perfect. We write $X\defeq \Spec(B)$. It induces two morphisms
\[[X/\calH]\overset{[f]}\rightarrow [X/\calG]\rightarrow [X/\calF]\]
of prestacks, that have perfect cotangent complexes over $[X/\calF]$. It induces two morphisms
\[[X/\calH]\overset{[f]}\rightarrow [X/\calG]\rightarrow [X/\calF]\]
of prestacks. The following is a generalization of \cite[3.50.]{lagthick}.
\begin{thm}\label{thm:lagfol}
	There is an equivalence:
	\[\Lag([f]/[X/\calF],n) \simeq \Lag(f/\calF,k[n](-2))\]
\end{thm}
\begin{proof}
	We denote $\calC \defeq \Mod_k^{\cfil}$. Let $C$ be the cofiber of
	\[\DR^\int(\calG/\calF)\otimes k(-2) \rightarrow \DR^\int(\calG/\calF).\]
	We claim that $|C|^\int$ is concentrated in external weights $0$ and $1$. The functor $\gr$ commutes with $|\vari|^\int$, so we have a fiber sequence:
	\[|\Sym_\calG(\cot^\int_{\calG/\calF}[-1]\langle 1 \rangle)\otimes k(-2)|^\int \rightarrow |\Sym_\calG(\cot^\int_{\calG/\calF}[-1]\langle 1 \rangle)|^\int \rightarrow \gr|C|^\int.\]
	But, by \cref{lm:relcot}, $\cot^\int_{\calG/\calF}$ is in filtration degrees greater than or equal to $1$, so $\Sym^m_\calG(\cot^\int_{\calG/\calF}[-1])$ is in filtration degrees greater than or equal to $m$. In that case, if we put $M_m \defeq \Sym^m_\calG(\cot^\int_{\calG/\calF}[-1])$, we observe that the filtration looks like
	\[\dots \rightarrow F^3 M_m \rightarrow F^2 M_m \overset\sim\rightarrow F^1 M_m \overset\sim\rightarrow F^0 M_m\]
	when $m \geq 2$. Then the filtration of $M_m\otimes k(-2)$ looks like
	\[\dots \rightarrow F^1(M_m\otimes k(-2)) = F^3M_m \rightarrow F^0(M_m\otimes k(-2)) = F^2 M_m.\]
	In both cases, the functor $|\vari|$ gives $F^2M_m$. So we have
	\[\gr^m|C|^\int \simeq \cofib(F^2 M_m \rightarrow F^2 M_m) \simeq 0.\]
	We then have:
	\[\Map_{\calC^{\cfil,\cfil}}(k\langle 2 \rangle[-n-2](0),C)\simeq 0\]
	hence
	\[\Map_{\calC^{\cfil,\cfil}}(k\langle 2 \rangle[-n-2](0),\DR^\int(\calG/\calF)\otimes k(-2)) \simeq \Map_{\calC^{\cfil,\cfil}}(k\langle 2 \rangle[-n-2](0),\DR^\int(\calG/\calF)).\]
	This also works for $\calH/\calF$.
	We have:
	\begin{align*}
		\Isot(f/\calF&,k[n](-2)) \\
			     &\simeq \Map_{\calC^{\cfil,\cfil}}(k\langle 2 \rangle[-n-2](0),\DR^\int(\calH/\calF)\otimes k(-2))\underset{\Map_{\calC^{\cfil,\cfil}}(k\langle 2 \rangle[-n-2](0),\DR^\int(\calG/\calF)\otimes k(-2))}\times *\\
			     &\simeq \Map_{\calC^{\cfil,\cfil}}(k\langle 2 \rangle[-n-2](0),\DR^\int(\calH/\calF))\underset{\Map_{\calC^{\cfil,\cfil}}(k\langle 2 \rangle[-n-2](0),\DR^\int(\calG/\calF))}\times *\\
			     &\simeq \Map_{\calC^{\cfil}}(k\langle 2 \rangle[-n-2],|\DR^\int(\calH/\calF)|^\int)\underset{\Map_{\calC^{\cfil}}(k\langle 2 \rangle[-n-2],|\DR^\int(\calG/\calF)|^\int)}\times * \\
	\end{align*}
	By \cref{thm:formformal}, this last line is equivalent to
	\[ \Map_{\calC^{\cfil}}(k\langle 2 \rangle[-n-2],\Gamma\left([X/\calF],\DR([X/\calH]/[X/\calF])\right))\underset{\Map_{\calC^{\cfil}}(k\langle 2 \rangle[-n-2](0),\Gamma\left([X/\calF],\DR([X/\calG]/[X/\calF])\right))}\times *\]
	which is, by definition, $\Isot([f]/[X/\calF],n)$. Now, assume $f$ has an isotropic structure $(\omega,h)$. We denote by $(\omega',h')$ the associated structure on $[f]$. We have a diagram:
	\[\begin{tikzcd}
		f^*\tan^\int_{\calG/\calF}[-n](2) \ar[r] \ar[d,"f^*\Theta_\omega"] & \tan^\int_{\calH/\calG}[-n+1](2) \ar[r]\ar[d,"\Theta_h"] & \tan^\int_{\calH/\calF}[-n+1](2)\ar[d,"(-1)^{-n+1}\Theta_h^\vee{[-n+1]}(2)"]\\
		f^*\cot^\int_{\calG/\calF} \ar[r] & \cot^\int_{\calH/\calF}\ar[r] & \cot^\int_{\calH/\calG}.
	\end{tikzcd}\]
	We can apply $B(0)\underset\calH\otimes\vari$ to this diagram and realize it. We get the diagram:
	\[\begin{tikzcd}
		\pi_\calG^*\tan_{[X/\calG]/[X/\calF]}[-n] \ar[r] \ar[d,"\pi_\calH^*\Theta_{\omega'}"] & \pi_\calH^*\tan_{[X/\calH]/[X/\calG]}[-n+1] \ar[r]\ar[d,"\Theta_{h'}"] & \pi^*\tan_{[X/\calH]/[X/\calF]}[-n+1]\ar[d,"(-1)^{-n+1}\Theta_{h'}^\vee{[-n+1]}(2)"]\\
		\pi_\calG^*\cot_{[X/\calG]/[X/\calF]} \ar[r] & \pi_\calH^*\cot_{[X/\calH]/[X/\calF]}\ar[r] & \pi_\calH^*\cot_{[X/\calH]/[X/\calG]}
	\end{tikzcd}\]
	where $\pi_\calG$ is the map $X\rightarrow [X/\calG]$ and $\pi_\calH$ is the map $X\rightarrow[X/\calH]$. We conclude because $\pi_\calH^*$ and $B(0)\underset\calH\otimes\vari$ reflect equivalences between perfect objects. So $(\omega,h)$ is non-degenerate if and only if $(\omega',h')$ is non-degenerate. We have:
	\[\Lag([f]/[X/\calF],n) \simeq \Lag(f/\calF,k[n](-2)).\]
\end{proof}
This implies the following:
\begin{prop}\label{prop:lagformdesc}
	If $X$ is a formal prestack, then
	\[\Symp(X,n) \simeq \lim\limits_{\Spec(A) \rightarrow X_\dR} \Symp(\calB_X(A)/\bD_{X_\dR}(A),k[n](-2)).\]
	Moreover, if $f : C \rightarrow X$ is a morphism of formal prestacks, there is a Cartesian square
	\[\begin{tikzcd}
		\Lag(f,n) \ar[r] \ar[d] & \Symp(X,n)\ar[d]\\
		\lim\limits_{\Spec(A)\rightarrow C_\dR}\Lag(f_\calB(A)^*/\bD_{C_\dR}(A),k[n](-2)) \ar[r] & \lim\limits_{\Spec(A)\rightarrow C_\dR}\Symp(f^*\calB_X(A)/\bD_{C_\dR}(A),k[n](-2))
		\arrow["\lrcorner"{anchor=center, pos=0.125}, draw=none, from=1-1, to=2-2]
	\end{tikzcd}\]
	where $f^* : \dPSt_{/X_\dR} \rightarrow \dPSt_{/C_\dR}$ is the restriction functor through $f$ and $f_\calB^* : f^*\calB_X \rightarrow \calB_C$ is the natural morphism obtained from $f$ on principal parts.
\end{prop}
\begin{proof}
	For the first equivalence, we know that $\DR(X) \simeq \lim\limits_{\Spec(A)\rightarrow X_\dR} \DR(X_A/\Spec(A))$ and non-degeneracy of a $2$-form on $X$ can be checked locally for each $X_A$, so
	\[\Symp(X,n) \simeq \lim\limits_{\Spec(A)\rightarrow X_\dR} \Symp(X_A/\Spec(A),n)\]
	Because $X$ is locally of finite presentation, we can restrict this limit to $\Spec(A) \rightarrow X_\dR$ such that $A$ is of finite presentation. By \cref{thm:lagfol} applied to
	\[\calF \defeq \bD_{X_\dR}(A) \rightarrow \calG \defeq \bD_{X_\dR}(A) \overset{f_A}\rightarrow \calH \defeq \calB_X(A),\]
	with $\pi_A : X_A \rightarrow \Spec(A) \simeq [f_A]$, we get:
	\begin{align*}
		\Symp(X_A/\Spec(A),n) &\simeq \Lag(\pi_A/\Spec(A),n+1) && \text{(see \cite[Example 2.26.]{ssgex})}\\
				      &\simeq \Lag(f_A/\bD_{X_\dR}(A),n+1) && \text{by \cref{thm:lagfol},}\\
				      &\simeq \Symp(\calB_X(A)/\bD_{X_\dR}(A),n) && \text{by \cite{ssgex} again.}\\
	\end{align*}
	For the second identity, $\Isot(f,n)$ fits into the pullback diagram
	\[\begin{tikzcd}
		\Isot(f,n) \ar[r] \ar[d] & \rm{pre}\Symp(X,n)\ar[d]\\
		* \ar[r,"0"] & \rm{pre}\Symp(C,n).
		\arrow["\lrcorner"{anchor=center, pos=0.125}, draw=none, from=1-1, to=2-2]
	\end{tikzcd}\]
	Restriction of the diagram $(\Spec(A) \rightarrow X_\dR)$ to the diagram $(\Spec(A) \rightarrow C_\dR)$ induces a commutative square:
	\[\begin{tikzcd}
		\rm{pre}\Symp(X,n) \ar[r] \ar[d] & \lim\limits_{\Spec(A)\rightarrow C_\dR} \rm{pre}\Symp(X_A/\Spec(A),n)\ar[d]\\
		\rm{pre}\Symp(C,n) \ar[r,"\sim"] & \lim\limits_{\Spec(A)\rightarrow C_\dR} \rm{pre}\Symp(C_A/\Spec(A),n).
	\end{tikzcd}\]
	If we take the fiber at $0$ of the right vertical map, we exactly get
	\[\lim\limits_{\Spec(A)\rightarrow C_\dR}\Isot(f_A/\Spec(A),n),\]
	and there is a natural pullback diagram obtained by the pasting lemma:
	\[\begin{tikzcd}
		\Isot(f,n) \ar[r] \ar[d] & \rm{pre}\Symp(X,n)\ar[d]\\
		\lim\limits_{\Spec(A)\rightarrow C_\dR}\Isot(f_A/\Spec(A),n) \ar[r] & \lim\limits_{\Spec(A)\rightarrow C_\dR} \rm{pre}\Symp(X_A/\Spec(A),n).
		\arrow["\lrcorner"{anchor=center, pos=0.125}, draw=none, from=1-1, to=2-2]
	\end{tikzcd}\]
	Because non-degeneracy can be checked locally on the family $(X_A/\Spec(A))$ for the symplectic structures, and on $(C_A/\Spec(A))$ for the Lagrangian structure, we get a Cartesian diagram:
	\[\begin{tikzcd}
		\Lag(f,n) \ar[r] \ar[d] & \Symp(X,n)\ar[d]\\
		\lim\limits_{\Spec(A)\rightarrow C_\dR}\Lag(f_A/\Spec(A),n) \ar[r] & \lim\limits_{\Spec(A)\rightarrow C_\dR}\Symp(X_A/\Spec(A),n)
		\arrow["\lrcorner"{anchor=center, pos=0.125}, draw=none, from=1-1, to=2-2]
	\end{tikzcd}\]
	we conclude by applying \cref{thm:lagfol} to the bottom line.
\end{proof}
\subsection{Descent for Lagrangian thickenings}
\par Because $\DR(X) \simeq \DR(X/X_\dR)$, $\DR$ satisfies descent with respects to the family $(X_A\rightarrow \Spec(A))_{\Spec(A)\rightarrow X_\dR}$. Formal thickenings satisfy the same kind of descent:
\begin{thm}\label{thm:descthick}
	For any convergent prestack $X$, we have
	\[\Thick(X) \simeq \lim\limits_{\Spec(A) \rightarrow X_\dR} \Thick(X_A)_{/\Spec(A)}.\]
\end{thm}
\begin{proof}
	If $X \rightarrow Y$ is a formal thickening of $X$, we have that $X_{\dR} \simeq Y_{\dR}$, so $Y_A$ is well defined as a family over $X_{\dR}$. One may easily check that the natural morphism $X_A \rightarrow Y_A$ is a formal thickening. Thus, we get a compatible family of formal thickenings and a functor
	\[Y_{\vari} : \Thick(X) \rightarrow \lim\limits_{\Spec(A) \rightarrow X_\dR} \Thick(X_A)_{/\Spec(A)}\]
	Conversely, if $(X_A \rightarrow Y_A)$ is a compatible family of formal thickenings over $\dAff/X_\dR$, then we set $Y' \defeq \colim Y_A$ and $Y \defeq Y'^\conv$. We claim that it defines a functor
	\[\colim : \lim\limits_{\Spec(A) \rightarrow X_\dR} \Thick(X_A)_{/\Spec(A)} \rightarrow \Thick(X).\]
	The map $X \rightarrow Y'$ is lafp by descent for colimits and because of commutation of colimits. We have $Y_C \simeq Y'^\wedge_{X_C}$ by descent as well. The functor $(\vari)_\red$ commutes with colimits, so we have $Y'_\red \simeq X_\red$, then $Y'_\dR\simeq X_\dR$. The whole point is to prove that $Y'$ has a deformation theory. Assume we have a Cartesian square
	\[\begin{tikzcd}
		B'\ar[d]\ar[r] & B\ar[d]\\
		A'\ar[r] & A
		\arrow["\lrcorner"{anchor=center, pos=0.125}, draw=none, from=1-1, to=2-2]
	\end{tikzcd}\]
	of commutative algebras, with $A'\rightarrow A$ a nilpotent extension, so $B' \rightarrow B$ is also a nilpotent extension. We have to prove that there exists a unique lift
	\[\begin{tikzcd}
		\Spec(A) \ar[r]\ar[d] & \Spec(A')\ar[d]\ar[rdd,bend left]\\
		\Spec(B) \ar[r] \ar[drr,bend right]& \Spec(B') \ar[rd,dashed] \\
						   && Y'
	\end{tikzcd}\]
	Such a lift will always be over $Y'_\dR\simeq X_\dR$ because $\Spec(B)_\dR \simeq \Spec(B')_\dR$. But $X_\dR$ has a deformation theory by \cref{lm:drdefth}, we then have a unique lift $\Spec(B') \rightarrow X_\dR$. The space of lifts
	\[\Spec(B') \rightarrow Y'\]
	over $X_\dR$ is equivalent to the space of sections
	\[\Spec(B') \rightarrow Y'\underset{X_\dR}\times \Spec(B') \simeq Y_{B'}.\]
	compatible with the commutative square. It is contractible since $Y_{B'}$ has a deformation theory. Now, $Y$
	\begin{itemize}
		\item is such that $X \rightarrow Y$ is lafp, because $Y' \rightarrow Y$ is lafp (see \cite[Remark 3.20.]{lagthick}),
		\item is convergent by construction,
		\item satisifies the unique lifting property along square zero extensions because $(\vari)^\conv$ commutes with limits, and
		\item $Y_C \simeq Y^\wedge_{X_C}$ for the same reason.
	\end{itemize}
	so $Y \in \Thick(X)$ and $Y_{(\vari)} \circ \colim \simeq \id$. We have that $\colim \circ Y_{(\vari)} \simeq \id$ because $X \simeq \colim\limits_{\Spec(A) \rightarrow X_\dR} X_A \simeq (\colim\limits_{\Spec(A) \rightarrow X_\dR} X_A)^\conv$ because it is convergent. This proves the proposition.
\end{proof}
\par Because $\DR$ satisfies this descent, we deduce:
\begin{prop}\label{prop:lagthickdesc}
	There is an equivalence:
	\[\Lag\Thick(X,n) \simeq \lim\limits_{\Spec(A)\rightarrow X_\dR} \Lag\Thick(X_A/\Spec(A),n).\]
\end{prop}
\begin{proof}
	By \cref{thm:descthick} and because $\DR(X) \simeq \lim \DR(X_A/\Spec(A))$, we see that
	\[\Isot\Thick(X,n) \simeq \lim\limits_{\Spec(A)\rightarrow X_\dR} \Isot\Thick(X_A/\Spec(A),n).\]
	Non-degeneracy is a local condition so this equivalence works also for Lagrangian structures.
\end{proof}
\section{Shifted Poisson structures}
\subsection{Reminder in the affine case}
\par Remember from \cite{lagthick} that, for a stable symmetric monoidal category $\calC$ (coming from a model category), there is an operad $\bP^\gr_n$ on $\calC^{\cfil}$ classifying Poisson algebras with a bracket of weight $-1$ and degree $1-n$. We can define in the same way an operad $\bP_{k[n](i)}$ classifying Poisson algebras with bracket of weight $i$ and degree $-n$. For instance, $\bP^\gr_n\simeq \bP_{k[n-1](-1)}$. These operads also exist in graded object $\calC^\gr$. We can also consider the operads $\bP^\gr_{k[n](m)\langle l \rangle}$ that exist in doubly filtered objects, or graded filtered objects, similar operads also exist if we add more filtrations or gradings. Deformation theory of these operads goes exactly the same way as for $\bP^\gr_n$. For $\calL \defeq k[n](i)$ and $A \in \Alg_{\bP_{\calL}}(\calC^{\cfil})$, we can endow
\[\Pol(A,\calL) \defeq \Hom_A\left(\Sym_A(\cot_A\otimes \calL[1]\langle -1 \rangle),A\right)\]
with a graded $\bP_{\calL[1]}$-algebra structure. If $A\in \CAlg(\calC^{\cfil})$, we also set:
\[\Pois(A,\calL) \defeq \Alg_{\bP_{\calL}}(\calC^{\cfil})\underset{\CAlg(\calC^{\cfil})}\times\{A\}\]
to be the space of \defterm{$\calL$-twisted Poisson structures} on $A$. Using deformation theory, we can prove the following:
\begin{prop}\label{prop:poismc}
	There is an equivalence:
	\[\Pois(A,\calL) \simeq \Map_{\Alg_{\Lie\{k\langle 1 \rangle\}}(\calC^{\gr,\cfil})}(k\langle 2 \rangle[-1],\Pol(A,\calL)\otimes\calL[1]).\]
\end{prop}
From this perspective, given an $\calL$-twisted Poisson structure $\pi$ on $A$, and if $A$ is of finite presentation, we naturally get a map
\[\Theta_\pi : \cot_A \otimes \calL \rightarrow \tan_A.\]
The Poisson structure is said to be \defterm{non-degenerate} if this is an equivalence. We denote by $\Pois^\nd(A,\calL)$ the maximal subgroupoid of $\Pois(A,\calL)$ consisting of non-degenerate structures. Following \cite[4.22.]{MSII} and \cite[Appendix A]{lagthick}, we have
\begin{prop}\label{prop:poisnd}
	There is an equivalence:
	\[\Pois^\nd(A,\calL)\simeq \Symp(A,\calL).\]
\end{prop}
If $A \in \Alg_{\bP_\calL}(\calC^{\cfil})$, it is possible to define its center $\rZ(A)$ that will be a $\bP_{\calL[1]}$-algebra, we refer to \cite[Section 3]{MSI} and \cite[Section 4]{lagthick} for a rigorous definition of the center. It is always given with a map of commutative algebras $\rZ(A) \rightarrow A$. We also remark that, if $A \in \CAlg(\calC^{\cfil})$, the graded algebra of $\calL$-twisted polyvectors
\[\Pol(A,\calL)\]
can be identified with the center of $A$:
\[\rZ(A)\]
where $A$ is viewed as a $\bP_{\calL\langle -1\rangle}$-algebra in $\calC^{\cfil,\gr}$ with zero Poisson bracket.
\begin{remark}
	The definition of the center highly depends on categorical models of complete filtered complexes as well as the brace construction. The author believes that it is possible to write a universal property for the center that relates it to Lurie's definition of centers. In the case $\bP_{k[n](0)}$, this has already been done by Safronov in \cite{Padd}. In general, one has to reprove a version of Poisson additivity in this context, which we believe is a straightforward generalization of the quoted article.
\end{remark}
\begin{lm}\label{lm:centfol}
	Let $\calF \rightarrow \calG$ be a morphism of perfect derived foliations over $A\in\CAlg_k$; assume that $\calG$ has a $k[n](-2)$-twisted Poisson structure in $\Mod_\calF(\Mod_k^{\cfil})$, then its center that we will denote
	\[\rZ_\calF(\calG)\]
	is a derived foliation over $A$.
\end{lm}
\begin{proof}
	By the trick described in \cite[5.1.]{lagthick}, it is always possible to endow $\rZ_\calF(\calG)$ with an extra filtration (that we will call ``external" filtration). If we denote by $\gr^\ext$ the graded pieces of this filtration, we have that $\gr^\ext\rZ_\calF(\calG) \simeq \Hom_\calG(\Sym_\calG(\tan_{\calG/\calF}\otimes\calL\langle -1\rangle),\calG)$. We have:
	\[\gr^\ext(\gr\rZ_\calF(\calG)) \simeq \gr(\gr^\ext\rZ_\calF(\calG))\]
	Because graded pieces are defined as cofibers, and cofibers commute with cofibers. Now,
	\[\gr(\gr^\ext\rZ_\calF(\calG))\simeq\gr\Hom_\calG(\Sym_\calG(\cot^\int_{\calG/\calF}\otimes k[n](-2)\langle -1 \rangle),\calG).\]
	By \cref{lm:grhom}, this is equivalent to
	\[\Hom_{\gr\calG}(\Sym_{\gr\calG}(\cot_{\gr\calG/\gr\calF}\otimes k[n](-2)\langle -1 \rangle),\gr\calG).\]
	And by \cref{lm:relcot}, this is equivalent to
	\[\Hom_{\gr\calG}(\Sym_{\gr\calG}(\cofib(\cot_\calF\rightarrow\cot_\calG)[-1](1)\underset{A}\otimes\gr\calG \otimes k[n](-2)\langle -1 \rangle),\gr\calG).\]
	Finally, after simplifying, we get
	\[\Sym_A(\fib(\tan_\calG\rightarrow\tan_\calF)[1-n](1)\langle 1 \rangle).\]
	We can view $\gr^i\rZ_\calF(\calF)$ as a filtered object. Now by our computation we know that $\gr^\ext (\gr^i\rZ_\calF(\calG))$ is fully concentrated in degree $i$. So 
	\begin{align*}
		\gr^i \rZ_\calF(\calG) &\simeq\gr^{\ext,i}(\gr^i(\rZ_\calF(\calG)))\\
				       &\simeq \gr^i\Sym_A(\fib(\tan_\calG\rightarrow\tan_\calF)[1-n](1)).
	\end{align*}
		Then, $\gr^0\rZ_\calF(\calG)\simeq A$, and $\rZ_\calF(\calG)$ is a derived foliation.
\end{proof}
It is also possible to define an operad $\bP_\calL^\rightarrow$ classifying relative Poisson algebras. Roughly speaking, a $\bP_\calL^\rightarrow$-algebra is the data of:
\begin{itemize}
	\item A $\bP_{\calL[-1]}$-algebra $A$,
	\item A $\bP_\calL$-algebra $B$, and
	\item A morphism of $\bP_\calL$-algebras $B \rightarrow \rZ(A)$.
\end{itemize}
If $\calL = k[n](-1)$, we recover the operad $\bP^\gr_{[n+2,n+1]}$ of \cite{lagthick}. Precomposing with map of commutative algebras $\rZ(A) \rightarrow A$ provides a functor
\[\Alg_{\bP^\rightarrow_\calL}(\calC^{\cfil}) \rightarrow \CAlg(\calC^{\cfil})^{\Delta^1}.\]
Given a map of commutative algebras $f : B \rightarrow A$, we set
\[\Cois(f,\calL) \defeq \Alg_{\bP^\rightarrow_\calL}(\calC^{\cfil})\underset{\CAlg(\calC^{\cfil})^{\Delta^1}}\times \{f\}\]
to be the space of \defterm{$\calL$-twisted coisotropic structures} on $f$. As before, it is possible to define what are \defterm{non-degenerate} coisotropic structures (see \cite[4.1.]{MSII}), we denote by $\Cois^\nd(f,\calL)$ the maximal subgroupoid of $\Cois(f,\calL)$ spanned by non-degenerate structures.
\begin{prop}\label{prop:coisnd}
	There is an equivalence:
	\[\Cois^\nd(f,\calL) \simeq \Lag(f,\calL).\]
\end{prop}
\par We will need the following lemma:
\begin{lm}\label{lm:centernd}
	Let $\calF\rightarrow \calG$ be a morphism of perfect derived foliations over $A$ of finite presentation; assume $\calG$ as a $\bP_{k[n](-2)}$-algebra structure in $\Mod_\calF(\Mod_k^{\cfil})$. Then the natural coisotropic structure on
	\[\rZ_\calF(\calG)\rightarrow \calG\]
	is non-degenerate.
\end{lm}
\begin{proof}
	With the external filtration, the natural coisotropic structure lifts into a non-degenerate $k[n+1](-2)\langle -1 \rangle$-twisted coisotropic structure on $\rZ_\calF(\calG) \rightarrow \calG\langle 0 \rangle$ by \cite[5.4.]{lagthick}. Note that $\rZ_\calF(\calG)$ is complete and almost of finite presentation (when viewed with the external filtration as well), so by the same argument than \cite[3.45.]{lagthick}, we know that
	\[|\Pol^\ext(\rZ_\calF(\calG),k[n+1](-2)\langle -1 \rangle)|^{\ext,\gr} \simeq \Pol(\rZ_\calF(\calG),k[n+1](-2))\]
	where:
	\begin{itemize}
		\item the functor $\Pol^\ext$ corresponds to polyvector fields applied to $\rZ_\calF(\calG)$ when taking into account both external and internal filtrations,
		\item it is given with an extra grading called here the Hodge grading,
		\item and $|\vari|^{\ext,\gr}$ is the functor $|\vari|^\ext$ (= realization for the external filtration), applied on each Hodge weights.
	\end{itemize}
	We reproduce it here for completeness. By construction of the center, there is a natural map
	\[|\Pol^\ext(\rZ_\calF(\calG),k[n+1](-2)\langle -1 \rangle)|^{\ext,\gr} \rightarrow \Pol(\rZ_\calF(\calG),k[n+1](-2))\]
	compatible with underlying $\Lie$ structures. On the $m$th graded piece we have:
	\[|(\tan_{\rZ_\calF(\calG)}^\ext[-n-2](2)\langle 1 \rangle)^{\otimes m}_{h\Sigma_m}|^\ext \rightarrow (\tan_{\rZ_\calF(\calG)}[-n-2](2))^{\otimes m}_{h \Sigma_m}\]
	where $\tan^\ext_{\rZ_\calF(\calG)}$ is the tangent computed by taking into account the external filtration on $\rZ_\calF(\calG)$. Because $\rZ_\calF(\calG)$ is almost of finite presentation, we know that $\tan^{\ext,\otimes m}$ can be computed in filtered complexes instead of complete filtered complexes (see the proof of \cite[3.45.]{lagthick}). Because $\tan_{\rZ_\calF(\calF)}^\ext\langle 1 \rangle$ is in external filtration degrees greater than or equal to $0$, the functor $|\vari|$ is the underlying object functor $(\vari)^u$. Taking cotangent in filtered objects commutes with $(\vari)^u$, so because $(\vari)^u$ is symmetric monoidal and $\cot^\ext_{\rZ_\calF(\calG)}$ is perfect, it also commutes with taking the tangent. We have
	\[|\tan_{\rZ_\calF(\calG)}^{\ext,\otimes m}|^\ext \simeq \tan_{|\rZ_\calF(\calG)|^\ext}^{\otimes m} \simeq \tan_{\rZ_\calF(\calG)}^{\otimes m}.\]
	We can also check that
	\[|\Pol^\ext(\calG\langle 0 \rangle,k[n+1](-2)\langle-1\rangle)|^\ext \simeq \Pol(\calG,k[n+1](-2)).\]
	We can see the same way as in the proof of \cite[3.50.]{lagthick} that non-degenerate structures are sent to non-degenerate structures. We conclude that the natural $k[n+1](-2)$-twisted coisotropic structure on
	\[\rZ_\calF(\calG)\rightarrow \calG\]
	is non-degenerate.
\end{proof}
\subsection{Global case}
\par Let $X$ be a formal prestack. We will define what are shifted Poisson structures on it. As discussed in \cref{sec:formloc}, we naturally have two presheaves on $X_\dR$: $\calB_X$ and $\bD_{X_\dR}$. We set:
\[\Pol(\calB_X/\bD_{X_\dR},k[n](-2)) : \begin{array}{rcl} \dAff^\op_{X_\dR} & \longrightarrow & \Alg_{\bP_{k[n+1](-2)\langle -1 \rangle}}(\Mod_k^{\gr,\cfil})\\
\Spec(A) & \longmapsto & \Pol(\calB_X(A)/\bD_{X_\dR}(A),k[n](-2)),\end{array}\]
and:
\[\Pol(X,n)\defeq \Gamma(X_\dR,\Pol(\calB_X/\bD_{X_\dR},k[n](-2))) \in \Alg_{\bP_{k[n+1](-2)\langle -1 \rangle}}(\Mod_k^{\gr,\cfil}).\]
\begin{defn}
	A $n$-shifted Poisson structure on $X$ is the data of a Maurer--Cartan element of $\Pol(X,n)$ of weight $\geq 2$:
	\[\Pois(X,n) \defeq \Map_{\Alg_{\Lie\{k\langle 1 \rangle\}}(\Mod_k^{\gr,\cfil})}\left(k\langle 2\rangle[-1],\Pol(X,n)\otimes k(-2)[n+1] \right).\]
	We say that it is \defterm{non-degenerate} if for any $\Spec(A) \rightarrow X_\dR$, $A$ of finite presentation, the $k[n](-2)$-twisted Poisson structure on $\calB_X(A)$ is non-degenerate, we denote by $\Pois^\nd(X,n)$ the maximal subgroupoid of $\Pois(X,n)$ spanned by non-degenerate $n$-shifted Poisson structures on $X$.
\end{defn}
By \cref{prop:poismc}, a shifted Poisson structure on a prestack $X$ is a compatible family of $\bD_{X_\dR}$-linear $\bP_{n+1}$-structures of filtration degree $-2$ on $\calB_X$. We can prove that we recover the correct definition for Artin stacks:
\begin{prop}
	If $X$ is an Artin stack locally of finite presentation, then
	\[\Pois(X,n) \simeq \Pois^{CPTVV}(X,n).\]
	The space of $n$-shifted Poisson structures on $X$ is the same as the space defined in \cite{CPTVV}.
\end{prop}
\begin{proof}
	It is because
	\begin{align*}
		\left|\Pol^{\geq 2}(\calB_X/\bD_{X_\dR},k[n](-2))\otimes k(-2)\right|^\int &\simeq \left|\Sym^{\geq 2}(\tan^\int_{\calB_X/\bD_{X_\dR}}[-n-1](2)\langle 1 \rangle)\otimes k(-2)\right|^\int\\
											   &\simeq \left(\Sym^{\geq 2}(\tan^\int_{\calB_X/\bD_{X_\dR}}[-n-1](2)\langle 1 \rangle)\otimes k(-2)\right)^{u,\int}\\
											   &\simeq \left(\Sym^{\geq 2}(\tan^\int_{\calB_X/\bD_{X_\dR}}[-n-1]\langle 1 \rangle)\right)^{u,\int}\\
											   &\simeq \left(\Pol^{\geq 2}(\calB_X/\bD_{X_\dR},k[n](0))\right)^{u,\int}.
	\end{align*}
	The first identity is by definition, the second one is because $\tan^\int_{\calB_X/\bD_{X_\dR}}$ is in internal weight $\geq -1$ (\cref{lm:relcot}), and the third one is because $(\vari)^u$ forgets filtration shifts. The functor $(\vari)^u$ is the same as the Tate realization functor of \cite{CPTVV}.
\end{proof}
\par If $f : C \rightarrow X$ is a morphism of prestacks, \cite{MSII} explains how to define coisotropic structures on $f$.
\begin{defn}
	The space $\Cois(X,n)$ of \defterm{$n$-shifted coisotropic structures} on $f$ is defined by the fiber product:
	\[\begin{tikzcd}
		\Cois(X,n) \ar[r] \ar[d] & \Pois(X,n)\ar[d]\\
		\Cois(f_\calB^*/\bD_{C_\dR},k[n](-2)) \ar[r] & \Pois(f^*\calB_X/\bD_{C_\dR},k[n](-2))
		\arrow["\lrcorner"{anchor=center, pos=0.125}, draw=none, from=1-1, to=2-2]
	\end{tikzcd}\]
	where $f^*_\calB$ is the natural morphism $f^*_\calB : f^*\calB_X \rightarrow \calB_C$. By definition, we also have a forgetful map
	\[\Cois(f,n) \rightarrow \Pois(C,n-1).\]
	We can define the space $\Cois^\nd(f,n)$ of \defterm{non-degenerate} $n$-shifted coisotropic structures by the fiber product
	\[\begin{tikzcd}
		\Cois^\nd(X,n) \ar[r] \ar[d] & \Pois^\nd(X,n)\ar[d]\\
		\Cois^\nd(f_\calB^*/\bD_{C_\dR},k[n](-2)) \ar[r] & \Pois^\nd(f^*\calB_X/\bD_{C_\dR},k[n](-2)).
		\arrow["\lrcorner"{anchor=center, pos=0.125}, draw=none, from=1-1, to=2-2]
	\end{tikzcd}\]
\end{defn}
\begin{prop}\label{prop:globnd}
	We have that
	\[\Pois^\nd(X,n) \simeq \Symp(X,n),\]
	and
	\[\Cois^\nd(f,n) \simeq \Lag(f,n).\]
\end{prop}
\begin{proof}
	We reduce to the affine case:
	\begin{align*}
		\Pois^\nd(X,n) &\simeq \lim\limits_{\Spec(A)\rightarrow X_\dR} \Pois^\nd(\calB_X(A)/\bD_{X_\dR}(A),k[n](-2)) && \text{by definition,}\\
			       &\simeq \lim\limits_{\Spec(A) \rightarrow X_\dR} \Symp(\calB_X(A)/\bD_{X_\dR}(A),k[n](-2)) && \text{by \cref{prop:poisnd},}\\
			       &\simeq \Symp(X,n) && \text{by \cref{prop:lagformdesc}.}
	\end{align*}
	The Cartesian diagram:
	\[\begin{tikzcd}
		\Lag(X,n) \ar[r] \ar[d] & \Symp(X,n)\ar[d]\\
		\Lag(f_\calB^*/\bD_{C_\dR},k[n](-2)) \ar[r] & \Symp(f^*\calB_X/\bD_{C_\dR},k[n](-2)).
		\arrow["\lrcorner"{anchor=center, pos=0.125}, draw=none, from=1-1, to=2-2]
	\end{tikzcd}\]
	and \cref{prop:coisnd} implies that
	\[\Cois^\nd(f,n) \simeq \Lag(f,n).\]
\end{proof}
\subsection{Two applications}
\subsubsection{Moduli of flat connections over non-proper varieties}
\par In \cite{TP2}, Pantev and To\"en build for $X$ a smooth non-proper complex variety a $(2-2d)$-shifted Poisson structure on the derived stack
\[\Vect^\nabla(X)\]
of flat connections on $X$ in a restricted sense. In fact, they build a morphism
\[f : \Vect^\nabla(X) \rightarrow \Vect^\nabla(\widehat\partial X)\]
with a $(3-2d)$-shifted isotropic structure. Fix $L/k$ a field extension, and an $L$-point $x : \Spec(L) \rightarrow \Vect^\nabla(X)$, consider the diagram:
\[\begin{tikzcd}
	* \ar[rrd, bend left]\ar[rd] \ar[ddr, bend right, "x"]\\
	&  \Vect^\nabla(X)^\wedge_x \ar[d] \ar[r,"f^\wedge_x"] & \Vect^\nabla(\widehat\partial X)^\wedge_x \ar[d] \\
	& \Vect^\nabla(X) \ar[r,"f"] & \Vect^\nabla(\widehat\partial X).
\end{tikzcd}\]
They prove that for any such $L$-point, the pullbacked $(3-2d)$-shifted isotropic structure on $f^\wedge_x$ is non-degenerate, that is, it defines a $(3-2d)$-shifted Lagrangian structure on $f^\wedge_x$. A corollary of \cref{prop:globnd} is that we get naturally a $(2-2d)$-shifted Poisson structure on each $\Vect^\nabla(X)^\wedge_x$, which was not known before. Indeed, in \cite{CPTVV}, they prove \cref{prop:poisnd} for $X$ an Artin derived stack, but here, $\Vect^\nabla(X)^\wedge_x$ is a formal neighborhood of a point, that is formal but not Artin. Although the cotangent complex of $\Vect^\nabla(X)$ is not perfect, we can still define what are Poisson structures on it; we don't know if it is possible to lift the structures on field points to a global Poisson structure. Anyway, To\"en suggested to the author that we could generalize his strategy but provided a working definition of a Tate derived stack with shifted symplectic and Poisson structures. We refer to \cite{tate} for discussions on this subject.
\subsubsection{BV Formalism}
\par In \cite[5.3.1.]{BVcois}, Grataloup gives an interpretation of BV reduction. He proves that Felder--Kazhdan's solution of the CME (\cite{CME}) defines, for $S : U \rightarrow \bA^1$ a function on a smooth affine scheme, a (strict) shifted Lagrangian correspondence
\[\begin{tikzcd}
	&\stCrit(S)\ar[rd]\ar[ld]\\
	\BV && \dCrit(S)
\end{tikzcd}\]
where $\stCrit(S)$ corresponds to the classical critical locus of $S$. Strictness is required here because he considers $\BV = \Spec(\CE(\calL))$, where $\CE(\calL)$ is the complete algebra of Chevalley--Eilenberg cochains of some $\calL_\infty$-algebroid and directly computes K\"ahler differentials on it, which is not a good model for the cotangent because this algebra is generally not cofibrant. We instead give another interpretation of the result of Felder and Kazhdan by directly building a derived foliation on $\stCrit(S)$ and prove that it induces a Lagrangian correspondence.
\begin{prop}\label{prop:BV}
	If $S : U \rightarrow \bA^1$ is a function on a smooth affine scheme $U=\Spec(A)$, and if $\stCrit(S)$ has a perfect cotangent complex, then there is a $(-1)$-shifted Lagrangian correspondence
	\[\begin{tikzcd}
	&\stCrit(S)\ar[rd]\ar[ld]\\
	\BV && \dCrit(S)
\end{tikzcd}\]
where $\BV$ is a formal thickening of $\stCrit(S)$. This induces a $(-1)$-shifted coisotropic structure on
\[\stCrit(S) \rightarrow \dCrit(S).\]
\end{prop}
Before proving this, we need a reminder on weak mixed graded modules.
\begin{defn}
	Recall that a weak mixed graded module $M$ is the data of a graded $k$-linear chain complex $M^i \in \rm{CH_k}$ with operations
	\[\eps_i : M^\bullet \rightarrow M^{\bullet+i}[1]\]
	satisfying $\eps_0 = d$ and
	\[\eps^2 = 0\]
	where $\eps \defeq \sum_i \eps_i$. Morphisms between two graded mixed modules $X$ and $Y$ are given by a list of arrows $f_i : X^\bullet \rightarrow Y^{\bullet+i}$ such that
	\[f \circ \eps_X = \eps_Y \circ f\]
	where $f \defeq \sum_i f_i$. We denote by $\rm{Ch}_k^{\eps,\gr}$ the $1$-category of weak mixed graded modules. It has a symmetric monoidal structure coming from the coalgebra structure of $k[\eps_i]/(\eps^2)$. We denote by $\CAlg^{\rm{st}}(\rm{Ch}_k^{\eps,\gr})$ the $1$-category of commutative algebras for this monoidal structure, where ``$\rm{st}$" stands for ``strict."
\end{defn}
Recall that by rectification for commutative algebras and by \cite[2.27.]{curvlie} (see also \cite{mixfil} for a model-independant proof), that, if we localize at gradedwise quasi-isomorphisms $W_{eq}$, we get an equivalence
\[\CAlg^{\rm{st}}(\rm{Ch}_k^{\eps,\gr})[W_{eq}^{-1}] \simeq \CAlg(\Mod_k^{\cfil})\]
compatible with the associated graded functors. We also need to prove the following lemma:
\begin{lm}\label{lm:drprod}
	Let $f : \Spec(B) \rightarrow X$ and $g : \Spec(B) \rightarrow Y$ be two morphisms of prestacks over $\Spec(B)_\dR$ having a cotangent complex, then
	\[\DR(\Spec(B)/(X\underset{\Spec(B)_\dR}\times Y)) \simeq \DR(\Spec(B)/X)\underset{\DR(B)}\otimes\DR(\Spec(B)/Y).\]
\end{lm}
\begin{proof}
	It suffices to check that the natural morphism
	\[\DR(\Spec(B)/X)\underset{\DR(B)}\otimes\DR(\Spec(B)/Y) \rightarrow \DR(\Spec(B)/(X\underset{\Spec(B)_\dR}\times Y))\]
	is an equivalence on associated graded. But on the left-hand side we have
	\begin{align*}\Sym(\cot_{\Spec(B)/X}[-1](1))\underset{\Sym(\cot_B[-1](1))}\otimes&\Sym(\cot_{\Spec(B)/Y}[-1](1))\\
	&\simeq \Sym(\cot_{\Spec(B)/X}[-1](1)\oplus_{\cot_B[-1](1)} \cot_{\Spec(B)/Y}[-1](1))\\
	&\simeq \Sym(\cot_{\Spec(B)/(X\underset{\Spec(B)_\dR}\times Y)}[-1](1))
	\end{align*}
	which is the right-hand side.
\end{proof}
The proof of \cref{prop:BV} is a bit complicated, especially because we have to properly describe the shifted Lagrangian structure, we only provide a sketch of the proof here and postpone a complete one for a future work.
\begin{proof}[Sketch of the proof of \cref{prop:BV}]
	Because $U$ is smooth, its cotangent complex is $\Omega^1_U$. We have an explicit model of $\dCrit(U)$:
	\[\dCrit(U) \simeq \Spec(\Sym_U(\Omega^1_U[1]),\iota_{\rm d_\dR S}).\]
	An explicit model for the $(-1)$-shifted symplectic form here is given by $\omega_U = \sum x_i\wedge p_i$ where $x_i$ are local coordinates and $p_i$ are phase coordinates. Let $\KT \defeq (\Sym_U(\calL[2]\oplus\Omega^1_U),\rm d = \iota_{\rm d_\dR S} + \rm d_{\rm{res}})$ be a Koszul--Tate resolution of $(\Sym_U(\Omega^1_U[1]),\rm d_\dR)$. This models the strict critical locus
	\[\stCrit(S) \simeq \Spec(\KT).\]
	As a graded module, $\calL$ is dualizable, so on the algebra
	\[\Sym_U(\calL[2]\oplus \calL^\vee[-1] \oplus \Omega^1_U)\]
	it is possible to define a $(-1)$-shifted Poisson bracket $\{\vari,\vari\}$ that is the sum of the Schouten--Nijenhuis bracket and the duality pairing of $\calL$ with $\calL^\vee$. If we complete at $\calL^\vee$, we get a filtered algebra
	\[\CE(\calL) \defeq \Sym_U(\calL[2]\oplus \widehat{\calL^\vee[-1]} \oplus \Omega^1_U).\]
	In \cite{CME}, Felder and Kazhdan define inductively a solution of what they call the ``Classical Master Equation." Explicitly, it is an element $\calS \in \CE(\calL)$ satisfying
	\begin{itemize}
		\item $\{\calS,\calS\} = 0$,
		\item $\calS_{|\calO_U} = S$, and
		\item $\calS = S + S_1 \text{ mod } \CE(\calL)^{\geq 2}$.
	\end{itemize}
	where $S_1$ is the composition
	\[k \rightarrow \calL\otimes \calL^\vee \overset{\rm d\otimes \id}\rightarrow \CE(\calL).\]
	Thus, $\{\calS,\vari\}$ defines a weak mixed graded structure on $\CE(\calL) \simeq \Sym_{\KT}(\KT\underset{A}\otimes \calL^\vee[-1](1))$ and, thanks to the reminder on weak mixed graded complexes, a derived foliation over $\stCrit(S)$. We write $\CE^\calS(\calL)$ this derived foliation. The duality bracket induces a ($\DR(\stCrit(S))$-linear) $k[-1](-2)$-shifted Lagrangian structure on
	\[\DR(\stCrit(S)) \rightarrow \CE(\calL)^\calS.\]
	Similarly, by \cref{thm:lagfol}, the $(-1)$-shifted symplectic structure on $\dCrit(S)$ gives a $k[-1](-2)$-twisted Lagrangian structure on
	\[\DR(\stCrit(S)) \rightarrow \DR(\stCrit(S)/\dCrit(S)).\]
	In fact, as in \cite[5.3.1.]{BVcois}, we get a $\DR(\stCrit(S))$-linear $k[-1](-2)$-twisted Lagrangian correspondence structure on
	\[\begin{tikzcd}
		& \KT\\
		\CE^\calS(\calL)\ar[ru] & & \DR(\stCrit(S)/\dCrit(S)).\ar[lu]
	\end{tikzcd}\]
	Because it is Lagrangian, and because we have assumed $\cot_{\stCrit(S)}$ to be perfect, $KT\underset{A}\otimes \calL^\vee$ with the differential is perfect, so $\CE(\calL)^\calS$ is a perfect derived foliation. So by \cref{lm:drprod} and \cref{thm:lagfol}, we get a $(-1)$-shifted Lagrangian correspondence (over $\stCrit(S)_\dR$)
	\[\begin{tikzcd}
		& \stCrit(S) \ar[rd]\ar[ld]\\
		\BV\defeq [\stCrit(S)/\CE^\calS(\calL)] && \dCrit(S).
	\end{tikzcd}\]	
	The $(-1)$-shifted coisotropic structure is because
	\[\stCrit(S) \rightarrow \overline\BV \times \dCrit(S)\]
	has a $(-1)$-shifted coisotropic structure by \cref{prop:globnd} (note that we couldn't apply \cite{CPTVV} here because $\BV$ is only formal), so we get a $(-1)$-shifted coisotropic correspondence
	\[\begin{tikzcd}
		& \stCrit(S) \ar[rd]\ar[ld]\\
		* && \overline\BV\times \dCrit(S).
	\end{tikzcd}\]
	but
	\[\begin{tikzcd}
		& \overline\BV\times\dCrit(S) \ar[rd]\ar[ld]\\
		\overline\BV\times\dCrit(S) && \dCrit(S)
	\end{tikzcd}\]
	is also a $(-1)$-shifted coisotropic correspondence; by composition of correspondences, we get a $(-1)$-shifted coisotropic structure on
	\[\stCrit(S) \rightarrow \dCrit(S).\]
\end{proof}
\par We plan to investigate BV formalism from this perspective in future work.
\section{AKSZ construction for shifted Poisson structures}
\subsection{Shifted Poisson structures are shifted Lagrangian thickenings}
\par In order to prove the AKSZ theorem for shifted Poisson structures, we need to prove the following theorem, which is a generalization of \cite[5.6]{lagthick}:
\begin{thm}\label{thm:lagthick}
	If $X$ is a formal prestack, there is an equivalence:
	\[\Pois(X,n) \simeq \Lag\Thick(X,n+1).\]
\end{thm}
Thus establishing a conjecture of Calaque (see \cite[Conjecture 3.4.]{cotstack})
\par To prove this theorem, we follow a relative analog of the ideas used in \cite{lagthick} to prove \cite[5.6.]{lagthick}. By definition, an $n$-shifted Poisson structure on $X$ is the data of a compatible family of $\bD_{X_\dR}$-linear $n$-shifted Poisson structures of weight $-2$ on $\calB_X$:
\[	\Pois(X,n) \simeq \lim\limits_{\Spec(A)\rightarrow X_\dR} \Pois(\calB_X(A)/\bD_{X_\dR}(A),k[n](-2)).\]
Thanks to \cref{prop:lagthickdesc}, we also know that
\[\Lag\Thick(X,n+1) \simeq \lim\limits_{\Spec(A)\rightarrow X_\dR} \Lag\Thick(X_A/\Spec(A),n+1).\]
But in our case ($X$ is formal), we can restrict this limit to $A$ of finite presentation. In that case, we have seen in \cref{sec:formloc} that $\calB(A)$ and $\bD_{X_\dR}(A)$ are both almost perfect derived foliations over $A^\red$ and that $X_A \simeq [\Spec(A^\red)/\calB_X(A)]$ and $\Spec(A) \simeq [\Spec(A^\red)/\bD_{X_\dR}(A)]$, to conclude, it suffices to prove the following:
\begin{prop}\label{prop:relthick}
	If $B \in \CAlg_k$ is of finite presentation, and $\calF \rightarrow \calG$ is a morphism in $\Fol_B$, then there is an equivalence
	\[\Pois(\calG/\calF,k[n](-2)) \simeq \Lag\Thick([\Spec(B)/\calG]/[\Spec(B)/\calF],n+1).\]
\end{prop}
\begin{proof}
	Like in \cite[5.6.]{lagthick}, we will need to introduce three intermediate spaces. Let
	\[\Cois(\calG/\calF,k[n](-2))\]
	be the spaces of diagrams
	\[\begin{tikzcd}
		\calH \ar[rr] && \calG\\
			      &\calF \ar[ur]\ar[ul]
	\end{tikzcd}\]
	in $\CAlg(\Mod_k^{\cfil})$ together with an $\calF$-linear $k[n](-2)$-twisted coisotropic structure on the map $\calH \rightarrow \calG$. Let
	\[\Cois^{nd,fol}(\calG/\calF,k[n](-2))\]
	be the spaces of diagrams
	\[\begin{tikzcd}
		\calH \ar[rr] && \calG\\
			      &\calF \ar[ur]\ar[ul]
	\end{tikzcd}\]
	in $\Fol_B$ together with an $\calF$-linear non-degenerate $k[n](-2)$-twisted coisotropic structure on the map $\calH \rightarrow \calG$. Finally, let
	\[\Lag^{fol}(\calG/\calF,k[n](-2))\]
	be the spaces of diagrams
	\[\begin{tikzcd}
		\calH \ar[rr] && \calG\\
			      &\calF \ar[ur]\ar[ul]
	\end{tikzcd}\]
	in $\Fol_B$ together with a $\calF$-linear $k[n](-2)$-twisted Lagrangian structure on the map $\calH \rightarrow \calG$. The last two spaces are equivalent by \cref{prop:coisnd}:
	\[\Lag^{fol}(\calG/\calF,k[n](-2)) \simeq \Cois^{nd,fol}(\calG/\calF,k[n](-2)).\]
	We can define, like in the proof of \cite[5.6.]{lagthick}, the ``center" functor:
	\[\rZ : \Alg_{\bP_{k[n](-2)}}(\Mod_\calF^{\cfil})^\fet \rightarrow \Alg_{\bP^\rightarrow_{k[n+1](-2)}}(\Mod_\calF^{\cfil})\]
	where
	\[\Alg_{\bP_{k[n](-2)}}(\Mod_\calF^{\cfil})^\fet\]
	denotes the subcategory of $\Alg_{\bP_{k[n](-2)}}(\Mod_\calF^{\cfil})$ with morphisms $B \rightarrow C$ such that the underlying morphism of commutative algebras is formally etale. This defines a map
	\[\Pois(\calG/\calF,k[n](-2)) \rightarrow \Cois(\calG/\calF,k[n+1](-2))\]
	that factors through $\Cois^{nd,fol}(\calG/\calG,k[n+1](-2))$ by \cref{lm:centfol} and \cref{lm:centernd}. We denote this factorization by $\rZ_\calF$. There is a map in the other direction
	\[\bP : \Cois^{nd,fol}(\calG/\calF,k[n+1](-2)) \rightarrow \Pois(\calG/\calF,k[n](-2))\]
	coming from the forgetful functor $\Cois(f,\calL) \rightarrow \Pois(C,\calL[-1])$ for an arbitrary map $f : D\rightarrow C$. There is an equivalence
	\[\bP\circ \rZ_\calF \simeq \id\]
	by construction of $\rZ_\calF$. Conversely, there is an equivalence
	\[\rZ_\calF \circ \bP \simeq \id.\]
	Indeed, if $C \rightarrow \calG$ has a non-degenerate $k[n+1](-2)$-twisted coisotropic structure, we naturally have a morphism of $\bP_{k[n+2](-2)}$-algebras
	\[C \rightarrow \rZ_\calF(\calG).\]
	It is an equivalence, because by non-degeneracy, $\cot_C \simeq \cot_{\rZ_\calF(\calG)}$. Hence we have:
	\begin{align*}
		\Pois(\calG/\calF,k[n](-2)) &\simeq \Cois^{nd,fol}(\calG/\calF,k[n+1](-2))\\
					    &\simeq \Lag^{fol}(\calG/\calF,k[n](-2))\\
					    &\simeq \Lag\Thick([\Spec(B)/\calG]/[\Spec(B)/\calF],n+1)&& \text{by \cref{thm:lagfol}.}
	\end{align*}
\end{proof}

\begin{proof}[Proof of \cref{thm:lagthick}]
	We have:
	\begin{align*}
		\Pois(X,n) &\simeq \lim\limits_{\Spec(A)\rightarrow X_\dR} \Pois(\calB_X(A)/\bD_{X_\dR}(A),k[n](-2))\\
			   &\simeq \lim\limits_{\Spec(A) \rightarrow X_\dR} \Lag\Thick(X_A/\Spec(A),n+1)\\
			   &\simeq \Lag\Thick(X,n+1).
	\end{align*}
	where the second line is \cref{prop:relthick} applied to $\calF = \bD_{X_\dR}(A)$ and $\calG = \calB_X(A)$.
\end{proof}
\subsection{The AKSZ construction}
\par We can finally state and prove the main theorem of this paper. Before, recall that a prestack $Y$ is said to be \defterm{$\calO$-compact} if for any $\Spec(A) \in \dAff$ the functor
\[\Gamma_A(Y\times \Spec(A),\vari) : \QCoh(Y\times\Spec(A)) \rightarrow \Mod_A\]
preserves colimits and perfect objects. It is said to be \defterm{$d$-oriented} if it is $\calO$-compact and equipped with a morphism
\[[Y] : \Gamma(Y,\calO_Y) \rightarrow k[-d]\]
such that for any $\calE$ perfect over $Y\times\Spec(A)$, the associated integration map:
\[\Gamma_A(Y\times\Spec(A),\calE) \rightarrow \Gamma_A(Y\times \Spec(A),\calE^\vee)^\vee[-d]\]
is an equivalence. Recall that, thanks to \cite[2.35. (2)]{sympgrpd}, if $L\rightarrow X$ has an $n$-shifted Lagrangian structure, and if $Y$ is $d$-oriented, then
\[\MapSt(Y,L) \rightarrow \MapSt(Y,X)\]
has an $(n-d)$-shifted Lagrangian structure.
\begin{thm}\label{thm:aksz}
	Let $X$ be an $n$-shifted derived Poisson formal prestack and $Y$ a $d$-oriented prestack. Assume, moreover, that $\MapSt(Y,X)$ is locally of finite presentation; then
	\[\MapSt(Y,X)\]
	has a natural $(n-d)$-shifted Poisson structure.
\end{thm}
\begin{proof}
	According to \cref{thm:lagthick}, we have an $(n+1)$-shifted Lagrangian thickening
	\[X\rightarrow X^\symp.\]
	After applying $\MapSt(Y,\vari)$ we get an $(n-d+1)$-shifted Lagrangian structure on
	\[f :\MapSt(Y,X) \rightarrow \MapSt(Y,X^\symp).\]
	It is a classical fact that $\calO$-compactness of $Y$ implies that both $\MapSt(Y,X)$ and $\MapSt(Y,X^\symp)$ have a perfect cotangent complex. By \cref{prop:mapdef} they both admit a deformation theory because $X$ and $L$ do. The map
	\[f : \MapSt(Y,X)\rightarrow \MapSt(X,X^\symp)\]
	is not lafp in general, but we can take the lafp approximation:
	\[f^\lafp : \MapSt(Y,X)^\lafp \rightarrow \MapSt(X,X^\symp)^\lafp.\]
	Because $\MapSt(Y,X)$ is assumed to be locally of finite presentation, hence lafp, we have that $\MapSt(Y,X) \simeq \MapSt(Y,X)^\lafp$. Moreover, $\MapSt(Y,X)$ is formal because it is locally of finite presentation. Note that $f^\lafp$ is itself lafp because a map between lafp prestacks is lafp. Now, we set
	\[\widehat{\MapSt(Y,X^\symp)} \defeq \MapSt(Y,X^\symp)^{\lafp,\wedge}_{\MapSt(Y,X)}\]
	be the formal completion of $\MapSt(Y,X)$ along $f^\lafp$. By $(1)$ of \cref{prop:lafp}, we know that $\MapSt(Y,X^\symp)^\lafp$ has a deformation theory. We also know that $f^\lafp$ is lafp. So the map $\MapSt(Y,X) \rightarrow \widehat{\MapSt(Y,X^\symp)}$ is a formal thickening.
	\par The prestack $\MapSt(Y,X^\symp)$ has a perfect cotangent complex, so $\MapSt(Y,X^\symp)^\lafp$ has a perfect cotangent complex. This implies that $f^\lafp$ has a perfect cotangent complex and by $(2)$ of \cref{prop:lafp} $f$ and $f^\lafp$ have the same relative cotangent complexes (because the source is lafp, so being equivalent on almost of finite presentation $k$-algebras is the same thing as being equivalent). We can fix:
	\[\MapSt(Y,X)^\symp \defeq \widehat{\MapSt(Y,X^\symp)}\]
	We get a formal thickening
	\[g : \MapSt(Y,X) \rightarrow \MapSt(Y,X)^\symp.\]
	Because $\MapSt(Y,X)^\symp \rightarrow \MapSt(Y,X^\symp)^\lafp$ is formally etale, $g$ has a relative cotangent complex equivalent to the one of $f$. Now $g$:
	\begin{itemize}
		\item has a natural $(n-d+1)$-shifted isotropic structure $(\omega,h)$ by transport given by the map
			\[\DR(\MapSt(Y,X^\symp)) \rightarrow \DR(\MapSt(Y,X^\symp)^\lafp)\rightarrow \DR(\MapSt(Y,X)^\symp),\]
		\item the structure is non-degenerate: we have that all vertical maps in the diagram
		\[\begin{tikzcd}
			g^*\tan_{\MapSt(Y,X)^\symp}[-n+d-1] \ar[r] \ar[d,"g^*\Theta_\omega"] & \tan_{\MapSt(Y,X)/\MapSt(Y,X)^\symp}[-n+d] \ar[r]\ar[d,"\Theta_h"] & \tan_{\MapSt(Y,X)}[-n+d]\ar[d]\\
			g^*\cot_{\MapSt(Y,X)^\symp} \ar[r] & \cot_{\MapSt(Y,X)}\ar[r] & \cot_{\MapSt(Y,X)/\MapSt(Y,X)^\symp}
		\end{tikzcd}\]
		are equivalences ($g$ and $f$ have the same relative cotangent complexes). So the structure is Lagrangian because $g^*$ reflects equivalences between perfect complexes (because $g$ is a formal thickening).
	\end{itemize}
	By the converse of $\cref{thm:lagthick}$, $\MapSt(Y,X)$ has an $(n-d)$-shifted Poisson structure.
\end{proof}
\begin{example}
	We can again use the same examples given in \cite[2.1.]{PTVV} for $d$-oriented prestack $Y$ such that $\MapSt(Y,X)$ is locally of finite presentation for $X$ formal:
	\begin{itemize}
		\item If $M$ is a compact connected oriented manifold of dimension $d$, then the Betti stack $Y \defeq M_B$ is $\calO$-compact and Poincaré duality gives a $d$-orientation on it. Compactness of $M$ implies that $\MapSt(Y,X)$ is locally of finite presentation.
		\item If $X$ is a smooth proper Deligne--Mumford stack of relative dimension $d$, with connected geometric fibers, the de Rham stack $X_\dR$ and Dolbeault stack $X_{Dol}\defeq B\widehat\tan_{X/k}$ are both $2d$-oriented. We can let $Y \defeq X_\dR$ or $Y \defeq X_{Dol}$ and we can reproduce the argument of \cite[19.1.3.1.(4)-(5)]{SAG} in order to see that $\MapSt(Y,X)$ is locally of finite presentation.
		\item If $Y$ is a Calabi--Yau smooth proper Deligne--Mumford stack of relative dimension $d$, with connected geometric fibers, then it is $d$-oriented. The prestack $\MapSt(Y,X)$ is also locally of finite presentation.
	\end{itemize}
\end{example}
\par We now give examples of new shifted Poisson structures.
\begin{example}
	When $G$ is an algebraic group over $k$ with associated Lie algebra $\fg$, after \cite[2.6.]{LiePois} and \cite[2.9.]{LiePois}, we know that:
	\begin{itemize}
		\item $\Pois(\rB G,2)\simeq \Sym^2(\fg)^G$, and
		\item $\Pois(\rB G, 1) \simeq \mathrm{QPois}(G)$ the $1$-groupoid of \textit{quasi-Poisson structures} on $G$.
	\end{itemize}
	So if $M$ is a compact connected manifold of dimension $d$, then the derived stack of $G$-local systems on $M$
	\[\mathrm{Loc}_G(M)\]
	carries a natural $(2-d)$-shifted Poisson structure, when $G$ is given with a $G$-equivariant $2$-tensor, and a $(1-d)$-shifted Poisson structure when $G$ is given with a quasi-Poisson structure.
\end{example}
\begin{remark}
	In the proof of \cref{thm:aksz}, we extract an $(n+1)$-shifted Lagrangian thickening of $X$ from an $n$-shifted Poisson structure on it, and we work with it. We believe that it is possible to work directly with derived foliations over $X$ instead. Indeed, an $n$-shifted Poisson structure over $X$ is the same thing as the data of a derived foliation $\calF$ over $X$ with a $k[n+1](-1)$-twisted Lagrangian structure on $\calF \rightarrow \calO_X$. As in the recent work of Alfieri (\cite{mapfol}), we expect that it is possible to define a derived foliation $\MapSt(Y,\calF)$ on
	\[\MapSt(Y,X)\]
	together with a $k[n-d+1](-1)$-twisted Lagrangian structure on
	\[\MapSt(Y,\calF) \rightarrow \calO_{\MapSt(Y,X)}.\]
	Hence an $(n-d)$-shifted Poisson structure on $\MapSt(Y,X)$.
\end{remark}
\printbibliography

\end{document}